\documentclass[EJP]{ejpecp} 
\DeclareMathOperator\erf{erf}%
\DeclareMathOperator\argmin{argmin}%

\SHORTTITLE{The Vervaat transform} 

\TITLE{The Vervaat transform of Brownian bridges and~Brownian~motion} 

\AUTHORS{%
 Titus~Lupu\footnote{Laboratoire de Math\'{e}matiques, Universit\'{e} Paris Sud, Orsay.
    \EMAIL{titus.lupu@u-psud.fr}}
  \and 
  Jim~Pitman\footnote{Statistics department, University of California, Berkeley.
  \EMAIL{pitman@stat.berkeley.edu}}
  \and
  Wenpin~Tang\footnote{Statistics department, University of California, Berkeley.
  \EMAIL{wenpintang@stat.berkeley.edu}}}

\KEYWORDS{Brownian quartet; Bessel bridges/processes; Markov property; path decomposition; semi-martingale property; size-biased sampling; Vervaat transform} 

\AMSSUBJ{60C05; 60J60; 60J65} 

\SUBMITTED{August 17, 2014}
\ACCEPTED{May 3, 2015}

\VOLUME{20}
\YEAR{2015}
\PAPERNUM{51}
\DOI{v20-3744}

\ABSTRACT{For a continuous function $f \in \mathcal{C}([0,1])$, define the Vervaat transform $V(f)(t):=f(\tau(f)+t \mod1)+f(1)1_{\{t+\tau(f) \geq 1\}}-f(\tau(f))$, where $\tau(f)$ corresponds to the first time at which the minimum of $f$ is attained. Motivated by recent study of quantile transforms of random walks and Brownian motion, we investigate the Vervaat transform of Brownian motion and Brownian bridges with arbitrary endpoints. When the two endpoints of the bridge are not the same, the Vervaat transform is not Markovian. We describe its distribution by path decomposition and study its semi-martingale property. The same study is done for the Vervaat transform of unconditioned Brownian motion, the expectation and variance of which are also derived.}

\begin{document}



\section{Introduction and main results}
In recent work of Assaf et al. \cite{AFP}, a novel path transform, called the quantile transform $Q$ has been studied both in discrete and continuous settings. 
Inspired by previous works in fluctuation theory, see e.g. Wendel \cite{Wendel} and Port \cite{Port}, the quantile transform of simple random walks is defined as follows. For $w$ a simple walk of length $n$ with increments of $\pm 1$,
the quantile transform associated to $w$ is defined by:
$$Q(w)_j:=\sum_{i=1}^j w(\phi_w(i))-w(\phi_w(i)-1) \quad \mbox{for}~1 \leq j \leq n,$$
where $\phi_w$ is the quantile permutation on $[1,n]$ 
defined by lexicographic ordering on pairs $(w(j-1),j-1)$, that is for $1 \leq i<j \leq n$, $w(\phi_w(i)-1)<w(\phi_w(j)-1)$, or $w(\phi_w(i)-1)=w(\phi_w(j)-1)$ and $\phi_w(i)<\phi_w(j)$. 

As shown in Assaf et al. \cite[Theorem $8.16$]{AFP}, the scaling limit of this transformation of simple random walks is the quantile transform of Brownian motion $B:= (B_t; 0 \leq t \leq 1)$:
\begin{equation}
\label{QBM}
Q(B)_t:=\frac{1}{2} L_1^{a(t)} + (a(t))^{+}-(a(t)-B_1)^{+} \quad \mbox{for}~0 \leq t\leq 1,
\end{equation}
where $(\cdot)^{+}:=\max(\cdot,0)$ is the positive part of any real number, $L_1^a$ is the local time of $B$ at level $a$ up to time $1$, and $a(t):=\inf\{a \in \mathbb{R}; \int_0^1 1_{B_{s} \leq a}ds>t\}$ is the quantile function of the occupation measure of Brownian motion. We refer readers to Dassios \cite{Dassios,Dassiosbis,Dassiostri}, Embrechts et al. \cite{ERY}, Yor \cite{Yorquant}, Bertoin et al. \cite{BCY}, Chaumont \cite{Chaumont1999} for further discussions on the quantile function $a(t)$ for $0 \leq t \leq 1$. See also the thesis of Forman \cite{Forman} for a detailed account of the quantile transform $Q$.

From the definition \eqref{QBM}, it is not quite clear what the quantile transform of Brownian motion looks like, and which properties it inherits from Brownian motion. 
The key result of  Assaf et al. \cite{AFP} is to identify the distribution of the somewhat mysterious $Q(w)$ with that of the Vervaat transform $V(w)$ defined as:
$$V(w)(i): =  \left\{ \begin{array}{ccl} w(\tau_V+i)-w(\tau_V)  &\mbox{for} & i \leq n-\tau_V, \\ w(\tau_V+i-n)+w(n)-w(\tau_V) &\mbox{for} & n-\tau_V \leq i \leq  n,  \end{array}\right.$$
where $\tau_V:= \min\{0 \leq j \leq n; w(j) \leq w(i)~\mbox{for}~0 \leq i \leq n \}$ is the first time at which the simple random walk reaches its global minimum. Consequently,
\begin{equation}
\label{identityQV}
(Q(B)_t; 0 \leq t \leq 1) \stackrel{(d)}{=} (V(B)_t; 0\leq t \leq 1)
\end{equation}
with
$$V(B)_t:=  \left\{ \begin{array}{ccl} B_{1-A+t}-B_{1-A}  &\mbox{for} & 0 \leq t \leq A, \\ B_{t-A}+B_{1}-B_{1-A} & \mbox{for} & A \leq t \leq 1,  \end{array}\right.$$
where $A$ is the almost sure arcsine split such that $1-A:=\argmin_{t \in [0,1]} B_t$. 
As a result of the identity \eqref{identityQV}, to understand the quantile transform of Brownian motion, it is equivalent to study its substitute, the Vervaat transform $V(B)$. 

Historically, Vervaat \cite{Vervaat} showed that if $B$ is conditioned to both start and
end at $0$, then $V(B)$ is normalized Brownian excursion:
\begin{theorem}\cite{Vervaat}
\label{Vervaat}
$(V(B^{0,br}); 0 \leq t \leq 1) \stackrel{d}{=} (B^{ex}; 0 \leq t \leq 1)$, where $(B^{0,br}_t; 0 \leq t \leq 1)$ is Brownian bridge of length $1$ starting at $0$ and ending at $0$, and $(B^{ex}_t; 0 \leq t \leq 1)$ is normalized Brownian excursion.
\end{theorem}
For usual definitions of Brownian bridge/excursion, we refer readers to Revuz and Yor \cite{RY}. 
Furthermore, Biane \cite{Biane} proved a converse theorem to Vervaat's result, that is to recover standard Brownian bridge from normalized Brownian excursion by uniform sampling:
\begin{theorem}\cite{Biane}
\label{Biane}
 Let $B^{ex}$ be normalized Brownian excursion, and $U$ be uniformly distributed random variable independent of $B^{ex}$. Then the shifted process $\theta(B^{ex},U)$ defined by
$$\theta(B^{ex},U)_t: =  \left\{ \begin{array}{ccl} B^{ex}_{U+t}-B^{ex}_U  &\mbox{for} & 0 \leq t \leq 1-U, \\ B^{ex}_{U+t-1}-B^{ex}_U &\mbox{for} & 1-U \leq t \leq 1,  \end{array}\right.$$
is standard Brownian bridge.
\end{theorem}
See also Pitman \cite{Pitman1999} as well as his monograph \cite[Exercise $6.1.1$ and $6.1.2$]{Pitman} for a simpler proof of these results using the cyclic lemma.

Chaumont \cite{Chaumont} extended Theorem \ref{Vervaat} to stable L\'evy processes. Chassaing and Jason \cite{CJ} considered a similar path transform of a reflected Brownian bridge conditioned on its local times at $0$. Miermont \cite{Miermont} proved a version of the theorem in the case of spectrally positive L\'evy processes, those are L\'evy processes with no negative jumps. Fourati \cite{Fourati} established the result for general L\'{e}vy processes under mild hypotheses, where a connection was made by Uribe \cite{Uribe2014} to L\'evy bridges conditioned to stay positive. Le Gall and Weill \cite{LeGall} studied a Vervaat-like transform of Brownian trees in terms of re-rooting. Theorem \ref{Vervaat} was extended to Markov processes by Fourati \cite{Fourati} and Vigon \cite{Vigon}, and to diffusion processes by Lupu \cite{Lupu}. 

By considering random times other than $\argmin$, Bertoin et al. \cite{BCP} studied the cyclic shift of Brownian bridges with non-zero endpoint. The more general case of  CEI processes, those are cyclically exchangeable increment processes, was treated by Chaumont and Uribe \cite{CU}. As far as we are aware, there has not been previous study of the Vervaat transform of unconditioned Brownian motion $B$, or of Brownian bridges $B^{\lambda,br}$ ending at $\lambda \ne 0$.

The contribution of the current paper is to provide path decomposition results of Vervaat transform of Brownian bridges, for simplicity, call them Vervaat bridges, with non-zero endpoint. Then we use the path decomposition to derive a collection of properties of the Vervaat transform of Brownian bridges and Brownian motion. The main result is stated as follows:
 \begin{theorem} 
\label{PT}
Let $\lambda<0$. Given $Z^{\lambda}$ the time of the first return to $0$ by $V(B^{\lambda,br})$, whose density is given by
\begin{equation}
\label{aldous}
f_{Z^{\lambda}}(t):=\frac{|\lambda|}{\sqrt{2 \pi t(1-t)^3}} \exp\left(-\frac{\lambda^2 t}{2(1-t)}\right),
\end{equation}
the path is decomposed into two (conditionally) independent pieces:
\begin{itemize}
\item
$(V(B^{\lambda,br})_{u}; 0 \leq u \leq Z^{\lambda})$ is Brownian excursion of length $Z^{\lambda}$;
\item
 $(V(B^{\lambda,br})_{u}; Z^{\lambda} \leq u \leq 1)$ is first passage bridge through level $\lambda$ of length $1-Z^{\lambda}$.
\end{itemize}
\end{theorem}
\begin{center}
\includegraphics[width=0.6\textwidth]{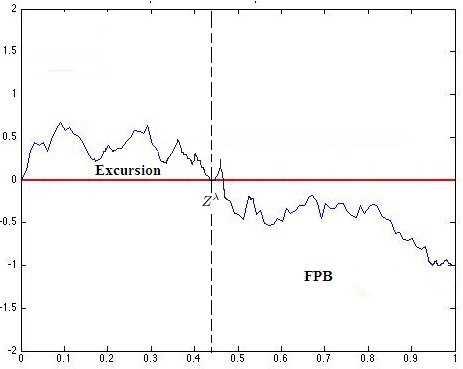}\\
Fig 1. Vervaat bridge $=$ Excursion $+$ First passage bridge.
\end{center}

Recall that Brownian excursion of length $l>0$ is defined by
$\left(B^{ex,l}_t:=\frac{1}{\sqrt{l}}B^{ex}_{t/l}; 0 \leq t \leq l \right)$, 
where $(B^{ex}_t;0 \leq t \leq 1)$ is normalized Brownian excursion. For $\lambda<0$, first passage bridge through level $\lambda$ of length $l>0$ is regarded as $(B_t; 0 \leq t \leq l)$ conditioned on $T_{\lambda}=l$, 
where $T_\lambda:=\inf\{t>0; B_t<\lambda\}$ is the first time at which Brownian motion hits $\lambda<0$.

Patrick Fitzsimmons points that Theorem \ref{PT} is closely related to Williams' decomposition of Brownian paths \cite{Williamsbis,Williams}. Actually, the result is a local version of Williams' decomposition, see e.g. Fitzsimmons \cite{Fitz} for further discussions. 
Note that $B^{\lambda,br}$ converges weakly to $B^{0,br}$ as $\lambda \rightarrow 0^{-}$ and the Vervaat transform $V$ is continuous at point $B^{0,br}$. Applying the continuous mapping theorem, Theorem \ref{Vervaat} is recovered as a weak limit $\lambda \rightarrow 0^{-}$ of Theorem \ref{PT}. 

The parametric family of densities $(f_{Z^{\lambda}})_{\lambda<0}$ appeared earlier in the work of Aldous and Pitman \cite[Corollary $5$]{AP}, when they studied the standard additive coalescent. For $B_1$ normally distributed with mean zero and unit variance,
\begin{equation}
\label{sbp}
Z^{\lambda} \stackrel{(d)}{=}\frac{B_1^2}{\lambda^2+B_1^2},
\end{equation} 
which corresponds to the mass of the size-biased tree component in the Brownian forest when the cutting intensity of the  Brownian CRT is $|\lambda|>0$. We refer readers to Pitman \cite[Chapter $4$ abd $10$]{Pitman} for further development.

For Vervaat bridges with positive endpoint $V(B^{\lambda,br})$ where $\lambda>0$, it is not hard to derive the following duality relation:
\begin{equation}
\left(V(B^{\lambda,br})_t; 0 \leq t \leq 1\right) \stackrel{(d)}{=} \left(V(B^{-\lambda,br})_{1-t}+\lambda; 0 \leq t \leq 1\right) \quad \mbox{for}~\lambda>0.
\end{equation}
In other words, looking backwards, we have the first piece of Brownian excursion above level $\lambda>0$, concatenated by the second one of first passage bridge from $\lambda$ to $0$. It is well-known that a first passage bridge from $\lambda>0$ to $0$ has the same distribution as the three dimensional Bessel bridge from $\lambda>0$ to $0$, see e.g. Biane and Yor \cite{BY}. We obtain the following decomposition of Vervaat bridges with positive endpoint:
\begin{corollary}
\label{PT2}
Let $\lambda>0$. Given $\widehat{Z}^{\lambda}$ the time of the last exit from $\lambda$ by $V(B^{\lambda,br})$, whose density is given by $f_{\widehat{Z}^{\lambda}}(t):=f_{Z^{-\lambda}}(1-t)\quad \mbox{as in \eqref{aldous}}$,
the path is decomposed into two (conditionally) independent pieces:
\begin{itemize}
\item
 $(V(B^{\lambda,br})_{u};0 \leq u \leq \widehat{Z}^{\lambda})$ is the three dimensional Bessel bridge of length $\widehat{Z}^{\lambda}$, starting from $0$ and ending at $\lambda$;
\item
$(V(B^{\lambda,br})_{u};\widehat{Z}^{\lambda} \leq u \leq 1)$ is Brownian excursion above level $\lambda$ of length $1-\widehat{Z}^{\lambda}$. 
\end{itemize}
\end{corollary}
\begin{center}
\includegraphics[width=0.6\textwidth]{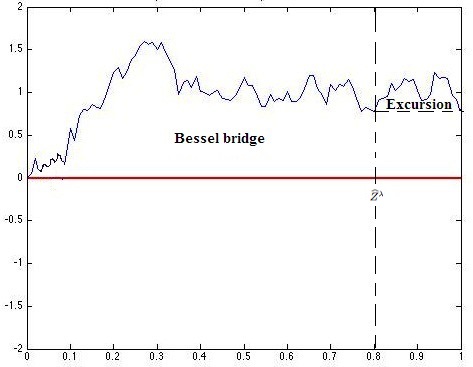}\\
Fig 2. Vervaat bridge $=$ Bessel bridge $+$ Excursion.
\end{center}
\textbf{Organization of the paper:} 
\begin{itemize}
\item
 In Section \ref{s2}, we provide two different proofs of Theorem \ref{PT}. The first proof is based on the weak limit approach,  where we make use of a bijection lemma proved by Assaf et al. \cite{AFP} (Subsection \ref{s21}).  The second proof relies on excursion theory, where the agreement formula by Pitman and Yor \cite{PY} is employed (Subsection \ref{s22}).
\item
In Section \ref{s3}, we give a thorough study of $V(B^{\lambda,br})$ where $\lambda \neq 0$ using Theorem \ref{PT} and Corollary \ref{PT2}. We prove that these processes are not Markov (Subsection \ref{s32}). However, they are semi-martingales and the canonical decompositions are given (Subsection \ref{s33}). Further we relate such processes to more elementary ones (Subsections \ref{s31} and \ref{s34}). The convex minorant of $V(B^{\lambda,br})$ is also studied (Subsection \ref{s35}). 
\item
In Section \ref{s4}, we focus on studying the Vervaat transform of Brownian motion. We prove that $V(B)$ is not Markov as well (Subsection \ref{s41}). We show that it is semi-martingale and the canonical decomposition is given (Subsections \ref{s42} and \ref{s43}). Finally, we provide explicit formulae for the first two moments of the Vervaat transform of Brownian motion (Subsection \ref{s44}).
\end{itemize}
\section{Path decomposition of Vervaat bridges}
\label{s2}
This section is devoted to proving Theorem \ref{PT}. First, we use a discrete approximation to obtain the path decomposition of $V(B^{\lambda,br})$ where $\lambda<0$. We obtain an analog to Theorem \ref{Biane} as a by-product. In the second part, we recover the same result via excursion theory. 
\subsection{Path decomposition via the weak limit approach}
\label{s21}
\subsubsection{Random walks analysis} 
\quad ~~We begin with the discrete time analysis of random walks, which is based on combinatorial principles. For a simple random walk $w$ of length $n$ with increments $\pm 1$, we would like to describe the law of $V(w^{a}):=(V(w)|w(n)=a)$, where $a<0$ has the same parity as $n$.

Recall that $\tau_V(w)=\min\{j \in [0,n]; w(j) \leq w(i)~\mbox{for}~ 0 \leq i \leq n \}$ is the first time at which the simple walk $w$ attains its minimum,  and $K(w)=n-\tau_V(w)$ is the distance from the position of the first minimum to the end. Following Assaf et al. \cite[Theorem $7.3$]{AFP}, the mapping $w \rightarrow (V(w), K(w))$ is a bijection between $walk(n)$, the set of simple walks of length $n$ and the set 
$$\{(v,k); v \in walk(n), v(j) \geq 0 ~\mbox{for} ~0 \leq j\leq k ~\mbox{and}~v(j)>v(n)~\mbox{for}~ k \leq j <n\},$$
where $k$, called a {\em helper variable}, records the splitting position in the original path. 
The following result is a direct consequence of this theorem related to Vervaat bridges.
\begin{lemma} \cite{AFP}
\label{FP}
Let $a<0$  have the same parity as $n$. Then $w^a \rightarrow (V(w^a), K(w^a))$ forms a bijection between $\{w \in walk(n): w(n)=a\}$ (simple bridges ending at $a<0$) and the set 
\begin{equation}
\label{setim}
\{(v,k); v \in walk(n), v(j) \geq 0 ~\mbox{for} ~0 \leq j\leq k ,~v(j)>a~\mbox{for}~ k \leq j <n ~\mbox{and}~v(n)=a\}.
\end{equation}
\end{lemma}
Observe that, to each pair $(v,k)$ in the set given by \eqref{setim}, one can associate a unique triple $(Z^{a}, f_{Z^{a}}^{br,1},f_{Z^{a}}^{br,2})$, where 
\begin{itemize}
\item $Z^{a}$ is the first time at which the path hits level $-1$;
\item  $f_{Z^{a}}^{br,1}$ is the path of first passage bridge of length $Z^{a}$ through level $-1$;
\item $f_{Z^{a}}^{br,2}$ is the path of first passage bridge of length $n-Z^{a}$ from $-1$ to $a<0$. 
\end{itemize}
Note that one may have the same triple $(Z^{a}, f_{Z^{a}}^{br,1},f_{Z^{a}}^{br,2})$ for different pairs $(v,k)$, since the mapping $w^a \rightarrow V(w^a)$ is not injective.

Now we compute the distribution of $Z^{a}$ by counting paths. The total number of simple bridges ending at $a<0$ is given by $\binom{n}{\frac{n+|a|}{2}}$. 
Fix $l>0$ an odd integer, 
the number of first passage bridges of length $l$ through level $-1$  is $\frac{1}{l} \binom{l}{\frac{l+1}{2}}$, and 
the number of first passage bridges of length $n-l$ starting from $-1$ to $a<0$ is $\frac{|a|-1}{n-l} \binom{n-l}{\frac{n-l+|a|-1}{2}}$,
see e.g. Feller \cite[Chapter III]{Feller}. Given $Z^a=l$, the total number of the Vervaat transform configurations is given by
$$\frac{|a|-1}{l(n-l)} \binom{l}{\frac{l+1}{2}} \binom{n-l}{\frac{n-l+|a|-1}{2}}.$$
By Lemma \ref{FP}, each Vervaat transform configuration $V(w^a)$ is counted $Z^a$ times. Hence,
\begin{equation}
\label{1}
\mathbb{P}(Z^{a}=l)=\frac{|a|-1}{n-l} \frac{\binom{l}{\frac{l+1}{2}} \binom{n-l}{\frac{n-l+|a|-1}{2}}}{\binom{n}{\frac{n+|a|}{2}}}.
\end{equation}

Combining the above arguments, we obtain the following path decomposition result for discrete Vervaat bridges with negative endpoint:
\begin{theorem} 
\label{discrete}
Let $a<0$ have the same parity as $n$. Given $Z^{a}:=\min\{j \geq 0;V(w^{a})_j=-1\}$, whose distribution is given by  \eqref{1}, the path is decomposed into two (conditionally) independent pieces:
\begin{itemize}
\item
$V(w^{a})|_{[0,Z^{a}]}$ is discrete first passage bridge of length $Z^{a}$ through level $-1$;
\item
$V(w^{a})|_{[Z^{a},n]}$ is discrete first passage bridge of length $n-Z^{a}$, starting from $-1$ and ending at $a$ .
\end{itemize}
\end{theorem}
\begin{center}
\includegraphics[width=0.6\textwidth]{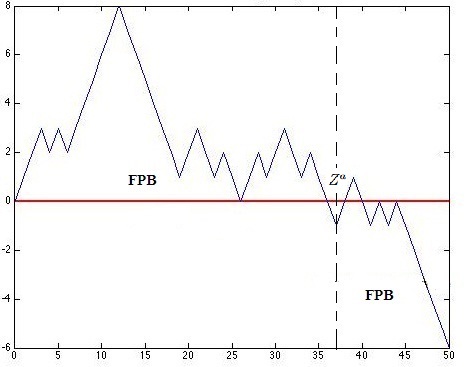}\\
Fig 3. Discrete Vervaat bridge $=$ First passage bridge $+$ First passage bridge.
\end{center}
\subsubsection{Passage to the weak limit}
\quad ~~We derive the path decomposition Theorem \ref{PT} from its discrete analog Theorem \ref{discrete} by appealing to invariance principles. 

For $\lambda<0$ and $0<t<1$, fix two sequences $(\lambda_n)_{n \in \mathbb{N}}$ and $(t_n)_{n \in \mathbb{N}}$ such that $\lambda_n \sim \lambda \sqrt{n}$ has the same parity as $n$, and $t_n := 2[\frac{tn}{2}]+1$. Let $S^{\lambda_n}$ be a simple random walk of length $n$ with increments $\pm 1$ which ends at $\lambda_n$, $V(S^{\lambda_n})$ be the associated discrete Vervaat bridge, and $Z^{\lambda_n}:= \inf \{j\geq0: V(S^{\lambda}_n)_j=-1\}$. 

Define $\left(S^{\lambda_n}(u); 0 \leq u \leq n \right)$ to be the linear interpolation of the discrete bridge $S^{\lambda_n}$, and $\left(V(S^{\lambda_n})(u); 0 \leq u \leq n \right)$ to be that of the discrete Vervaat bridge $V(S^{\lambda_n})$. In the next lemma, we recall some invariance principles in the metric space $\mathcal{C}[0,1]$, that is the space of continuous paths on $[0,1]$. For background on the weak convergence in $\mathcal{C}[0,1]$, we refer readers to Billingsley \cite[Chapter 2]{Bill2}.
\begin{lemma}
~
\begin{enumerate}
\item 
$(\frac{1}{\sqrt{n}}V(S^{\lambda_n})(nu); 0 \leq u \leq 1)$ converges weakly to $(V(B^{\lambda,br})_u; 0 \leq u \leq 1)$  in $\mathcal{C}[0,1]$.
\item
Given $Z^{\lambda_n}=t_n$, $(\frac{1}{\sqrt{n}}V(S^{\lambda_n})(nu); 0 \leq u \leq t)$ converges weakly to Brownian excursion of length $t$ in $\mathcal{C}[0,1]$, and $(\frac{1}{\sqrt{n}}V(S^{\lambda_n})(nu); t \leq u \leq 1)$ converges weakly to first passage bridge through level $ \lambda$ of length $1-t$ in $\mathcal{C}[0,1]$, (conditionally) independent of the excursion.
\end{enumerate}
\end{lemma}
\begin{proof} It is well-known that $(\frac{1}{\sqrt{n}}S^{\lambda_n}(nu);0 \leq u \leq 1)$ converges weakly to $(B^{\lambda,br}_u;0 \leq u \leq 1)$ in $\mathcal{C}[0,1]$. Note that the Vervaat transform $V: \mathcal{C}[0,1] \rightarrow \mathcal{C}[0,1]$ is continuous at point $B^{\lambda,br}$, which has a unique minimum on $[0,1]$ with probability $1$. The assertion $(1)$ follows readily from the continuous mapping theorem, see e.g. Billingsley \cite[Theorem $5.1$]{Bill2}. According to Theorem \ref{discrete}, given $Z^{\lambda_n}=t_n$, the path of $V(S^{\lambda_n})$ is split into two (conditionally) independent pieces of discrete first passage bridges. Following from  Iglehart \cite{Iglehart} and Bertoin et al. \cite{BCP}, the scaled first passage bridge through level $-1$ converges weakly to Brownian excursion in $\mathcal{C}[0,1]$, and the scaled first passage bridge from $-1$ to $\lambda_n$ converges weakly to first passage bridge through level $\lambda$ in $\mathcal{C}[0,1]$. This proves $(2)$.  
\end{proof}

From \eqref{1}, we obtain:
\begin{equation}
\label{3}
n \mathbb{P}(Z^{{\lambda}_n}=t_n)= \frac{n|\lambda_n|}{n-t_n} \frac{\binom{t_n}{\frac{t_n+1}{2}} \binom{n-t_n}{\frac{n-t_n+|\lambda_n|-1}{2}} }{\binom{n}{\frac{n+|\lambda_n|}{2}}}.
\end{equation}
Using Stirling's formula, we have that
$$ \binom{n}{\frac{n+|\lambda_n|}{2}} \sim \sqrt{\frac{2}{\pi n}} 2^n \exp \left(-\frac{\lambda^2}{2}\right);$$
$$\binom{t_n}{\frac{t_n+1}{2}} \sim \sqrt{\frac{2}{\pi nt}} 2^{nt};$$
and
$$\binom{n-t_n}{\frac{n-t_n+|\lambda_n|-1}{2}} \sim \sqrt{\frac{2}{\pi n(1-t)}} 2^{n(1-t)} \exp\left(-\frac{\lambda^2}{2(1-t)}\right).$$
Injecting these terms in \eqref{3}, we deduce the limiting distribution as $n \rightarrow \infty$ given by \eqref{aldous} . By a  local limit argument, see e.g. Billingsley \cite[Exercise $25.10$]{Bill}, we conclude that $Z^{\lambda}$ has density $f_{Z^{\lambda}}$ given as in \eqref{aldous}.

The next theorem is a direct consequence of Theorem \ref{PT} and should be called a corollary at best. Because of its importance, however, we give it status of a theorem.
\begin{theorem}
Let $\lambda<0$. Given $Z^{\lambda}$ the time of the first return to $0$ by $V(B^{\lambda,br})$, the split position $A^{\lambda}:=1-\argmin_{t \in [0,1]} B^{\lambda,br}_t$ is (conditionally) independent of $V(B^{\lambda,br})$, and is  uniformly distributed on $[0,Z^{\lambda}]$, In particular, its density is
$$f_{A^{\lambda}}(a):=\int_a^1 \frac{f_{Z^{\lambda}}(t)}{t}dt,$$
where $f_{Z^{\lambda}}$ is given by \eqref{aldous}.
\end{theorem}
\begin{proof} Given a Vervaat bridge configuration $V(w^{\lambda_n})$, the {\em helper variable} $K(w^{\lambda_n})$ takes values in $\{0,..., Z^{\lambda_n}-1\}$, where $Z^{\lambda_n}$ is the first time at which the path hits $-1$. This implies that given $Z^{\lambda_n}$, the distance from the minimum position of the original bridge to the end is uniformly distributed on $[0,Z^{\lambda_n}]$. The results are obtained by passing to the weak limit. 
\end{proof}
\begin{corollary}
\label{Bianebis}
Let $\lambda<0$. Let $Z^{\lambda}$ be the time of the first return to $0$ by $V(B^{\lambda,br})$ and $A^{\lambda}$ be uniformly distributed on $[0,Z^{\lambda}]$. Then the shifted process $\theta(V(B^{\lambda,br}),A^{\lambda})$ defined as in Theorem \ref{Biane} is Brownian bridge ending at $\lambda$, which attains its minimum at $1-A^{\lambda}$.
\end{corollary}
Again by the weak convergence of $B^{\lambda,br}$ to $B^{0,br}$ as $\lambda \rightarrow 0^{-}$ and the continuity of $V$ at point $B^{0,br}$, Theorem \ref{Biane} is recovered as a weak limit $\lambda \rightarrow 0^{-}$ of Corollary \ref{Bianebis}. The above corollary is extended to $\lambda \leq 0$.
\subsection{Path decomposition via excursion theory}
\label{s22}
In the current subsection, we provide an alternative proof of Theorem \ref{PT} using excursion theory. The proof relies on the decomposition of bridges at their minimum, a variant of the decomposition at the maximum that appeared in Pitman and Yor \cite{PY}. Fix $\lambda<0$, we begin with some notations.

For $w \in \mathcal{C}[0,\infty)$, that is the space of continuous paths on $[0,\infty)$, define the lifetime of the path $w$ by $\zeta(w):=\inf\{t \geq 0; w_t=\Delta\}$,
where $\Delta$ is a cemetery point. Let $$\mathcal{C}_f:=\{w \in \mathcal{C}[0,\infty); \zeta(w)<\infty\}$$ be the space of continuous paths of finite length. Given a distribution $\mathbb{Q}$ on $\mathcal{C}_f$, $\mathbb{Q}^{\wedge}$ is the image by time reversal: for $F$ Borel measurable,
$$\mathbb{Q}^{\wedge}[(w_s; 0 \leq s \leq \zeta(w)) \in F]:=\mathbb{Q}[(w_{\zeta(w)-s}; 0 \leq s \leq \zeta(w)) \in F].$$
Given $\mathbb{Q}$ and $\mathbb{Q}'$ two distributions on $\mathcal{C}_f$, $\mathbb{Q}\circ \mathbb{Q}'$ is the distribution obtained by concatenating two independent paths, one distributed as $\mathbb{Q}$ and the other as $\mathbb{Q}'$: for $F', F''$ Borel measurable and $F=F' \otimes F''$,
$$\mathbb{Q}\circ \mathbb{Q}'(w^1 \otimes w^2 \in F):=\mathbb{Q}(w^1 \in F') \times \mathbb{Q}(w^2 \in F''),$$
where $w^1 \otimes w^2(s):=w_s^1 1_{0 \leq s \leq \zeta(w^1)}+w_{s-\zeta(w^1)}^2 1_{\zeta(w^1) < s \leq \zeta(w^1)+\zeta(w^2)}$ for $0 \leq s \leq \zeta(w^1)+\zeta(w^2)$, is the concatenation of two paths.

Let $p_{t}(x,y)$ be the heat kernel defined as
$$p_{t}(x,y):=\frac{1}{\sqrt{2\pi t}}\exp\left(-\frac{(y-x)^{2}}{2t}\right),$$
and $\mathbb{P}^{T}_{0,\lambda}$ be the distribution of Brownian bridge of length $T$ from $0$ to $\lambda$, and $\mathbb{P}_{x}^{T_{y}}$ be the distribution of Brownian motion starting from $x$ until the first time at which it hits $y$ for $y<x$. 

The following agreement formula is read from Pitman and Yor \cite[Corollary $3$]{PY}:
\begin{equation}
\label{PitmanYor}
\int_{0}^{+\infty}dT\,p_{T}(0,\lambda)\mathbb{P}^{T}_{0,\lambda}=2\int_{-\infty}^{\lambda}
dy\,\mathbb{P}_{0}^{T_{y}}\circ\mathbb{P}_{\lambda}^{T_{y}\wedge},
\end{equation}
where the LHS of \eqref{PitmanYor} follows from L\'evy-It\^o's description of Brownian excursions, and the RHS of \eqref{PitmanYor} stems from Williams' decomposition of Brownian paths. We also refer readers to Yen and Yor \cite[Chapter $6$]{YY} for details.
Furthermore, $\mathbb{P}_{0}^{T_{y}}$ can be decomposed as:
\begin{displaymath}
\mathbb{P}_{0}^{T_{y}}=\mathbb{P}_{0}^{T_{y-\lambda}}\circ\mathbb{P}_{y-\lambda}^{T_{y}}.
\end{displaymath}
Therefore,
\begin{equation}
\label{Eq1}
\int_{0}^{+\infty}dT\,p_{T}(0,\lambda)\mathbb{P}^{T}_{0,\lambda}=2\int_{-\infty}^{\lambda}
dy\,\mathbb{P}_{0}^{T_{y-\lambda}}\circ
\mathbb{P}_{y-\lambda}^{T_{y}}\circ\mathbb{P}_{\lambda}^{T_{y}\wedge}.
\end{equation}

Now we extend the definition of the Vervaat transform to $\mathcal{C}_f$: given a continuous function $f$ on $[0,T]$ and $\tau(f)$ the first time at which it attains its minimum, define
\begin{displaymath}
V_T(f)(t):=f(\tau(f)+t \mod T)+f(T)1_{\lbrace t+\tau(f)\geq T\rbrace}-f(\tau(f)).
\end{displaymath}
Applying the Vervaat transform $V_T$ to \eqref{Eq1}, we obtain:
\begin{align}
\int_{0}^{+\infty}dT\,p_{T}(0,\lambda)V_{T}(\mathbb{P}^{T}_{0,\lambda})&=2\int_{-\infty}^{\lambda}
dy\,\mathbb{P}_{\lambda-y}^{T_{0}\wedge}\circ\mathbb{P}_{\lambda-y}^{T_{0}}\circ
\mathbb{P}_{0}^{T_{\lambda}}\notag\\
&\label{Eq2}=
2\left(\int_{0}^{+\infty}
dy\,\mathbb{P}_{y}^{T_{0}\wedge}\circ\mathbb{P}_{y}^{T_{0}}\right)\circ
\mathbb{P}_{0}^{T_{\lambda}}.
\end{align}
By taking $\lambda=0$ in \eqref{Eq2}, we see that
\begin{displaymath}
\int_{0}^{+\infty}dT\,p_{T}(0,0)V_{T}(\mathbb{P}^{T}_{0,0})=
2\int_{0}^{+\infty}
dy\,\mathbb{P}_{y}^{T_{0}\wedge}\circ\mathbb{P}_{y}^{T_{0}}.
\end{displaymath}

Let $\mathbb{Q}^{T}_{0,0}$ be the distribution of  Brownian excursion of length $T$, i.e. the three dimensional Bessel bridge from $0$ to $0$. A slight modification of Vervaat's result \cite{Vervaat} gives that $V_{T}(\mathbb{P}^{T}_{0,0})=\mathbb{Q}^{T}_{0,0}$. Thus,
\begin{equation}
\label{332}
2\int_{0}^{+\infty}
dy\,\mathbb{P}_{y}^{T_{0}\wedge}\circ\mathbb{P}_{y}^{T_{0}}=
\int_{0}^{+\infty}dT\,p_{T}(0,0)\mathbb{Q}^{T}_{0,0}.
\end{equation}
Injecting \eqref{332} into \eqref{Eq2}, we have that
\begin{equation}
\label{Eq3}
\int_{0}^{+\infty}dT\,p_{T}(0,\lambda)V_{T}(\mathbb{P}^{T}_{0,\lambda})=
\left(\int_{0}^{+\infty}
ds\,p_{s}(0,0)\mathbb{Q}^{s}_{0,0}\right)\circ
\mathbb{P}_{0}^{T_{\lambda}}.
\end{equation}
By disintegrating \eqref{Eq3} with respect to the lifetime of paths, we see that $V_{T}(\mathbb{P}^{T}_{0,\lambda})$ is a concatenation of Brownian excursion and Brownian first passage bridge.

Let $g_t(\lambda)$ be the density of the first hit to $\lambda$ by Brownian motion starting from $0$:
\begin{equation}
\label{Eq4}
g_t(\lambda):=\dfrac{\vert\lambda\vert}{\sqrt{2\pi t^{3}}}
\exp\left(-\dfrac{\lambda^{2}}{2t}\right).
\end{equation}
It follows from \eqref{Eq3} that the density of the splitting position $Z^{\lambda}$ between the excursion and the first passage bridge in $V(\mathbb{P}^{1}_{0,\lambda})$ is given by
\begin{displaymath}
\dfrac{p_{t}(0,0)g_{1-t}(\lambda)}{p_{1}(0,\lambda)}\,= f_{Z^{\lambda}}(t) \quad \mbox{defined as in}~\eqref{aldous}.
\end{displaymath}
\section{The Vervaat transform of Brownian bridges}
\label{s3}
\quad In this section, we study thoroughly Vervaat bridges with non-zero endpoint. First, we give an alternative construction of $V(B^{\lambda,br})$ where $\lambda \neq 0$ via Brownian bridges conditioned on its local times at $0$. Next we show that these processes are not Markovian. They are semi-martingales and canonical decompositions are given in both positive and negative endpoint cases. We further relate Vervaat bridges to drifting excursion by additional conditioning. To close the section, we study some properties of the convex minorant of $V(B^{\lambda,br})$ where $\lambda<0$.
\subsection{Construction of Vervaat bridges via Brownian bridge conditioned on its local times}
\label{s31}
In this subsection, we provide an alternative construction of Vervaat bridges with negative endpoint in terms of Brownian bridge which ends at $0$. It is straightforward that Vervaat bridges with positive endpoint can be treated similarly by time reversal.

Let $\lambda<0$. By Theorem \ref{PT}, conditioned on $Z^{\lambda}$ the time of the first return to $0$, $V(B^{\lambda,br})$ is split into Brownian excursion of length $Z^{\lambda}$, followed by first passage bridge of length $1-Z^{\lambda}$ through $\lambda$, independent of each other. Note that $V(B^{\lambda,br})$ looks like first passage bridge of length $1$ through $\lambda$, except that it starts with a piece of excursion.

Recall that first passage bridges of length $1$ can be constructed via standard Brownian bridge by conditioning on its local times. Let $(F^{\lambda}_t; 0 \leq t \leq 1)$ be a first passage bridge of length $1$ through $\lambda<0$. Following from Bertoin et al. \cite{BCP}, $F^{\lambda}$ has the same distribution as
\begin{equation}
\label{4}
(|B^{0,br}_t|-L_t^0(B^{0,br}); 0 \leq t \leq 1) \quad \mbox{conditioned on}~L_1^0(B^{0,br})=|\lambda|,
\end{equation}
where $L_t^0$ is the local times of Brownian bridge $B^{0,br}$ at level $0$ up to time $t$. We refer readers to Aldous \cite{Aldous}, and Chassaing and Janson \cite{CJ} for discussions of Brownian bridge conditioned on its local times.
In light of the above construction, we obtain:
\begin{theorem}
Let $X$ be the process whose distribution is that of $B^{0,br}$ conditioned on $L_1^0(B^{0,br})=|\lambda|$, $U$ be uniformly distributed on $(0,1)$ independent of $X$, and $(G_U,D_U)$ be the signed excursion interval which contains $U$. 
Let $\widetilde{X}$ be the process obtained by exchanging the position of the excursion of $X$ straddling time $U$ and the path along $[0,G_U]$, namely:
$$\widetilde{X}_t =     \left\{ \begin{array}{rcl}  X_{t+G_U} & \mbox{for} & 0 \leq t \leq D_U-G_U \\ X_{t-D_U+G_U}  & \mbox{for} & D_U-G_U \leq t \leq D_U \\  X_t & \mbox{for} &  D_U \leq t \leq 1.\end{array}\right.$$
Then for $\lambda<0$, $(V(B^{\lambda,br})_t; 0 \leq t \leq 1)$ has the same distribution as
\begin{equation}
\label{5}
(|\widetilde{X}(t)|-L_t^0(\widetilde{X}); 0 \leq t \leq 1) \quad \mbox{conditioned on}~L_1^0(\widetilde{X})=|\lambda|.
\end{equation}
\end{theorem}
\begin{center}
\includegraphics[width=1.1\textwidth]{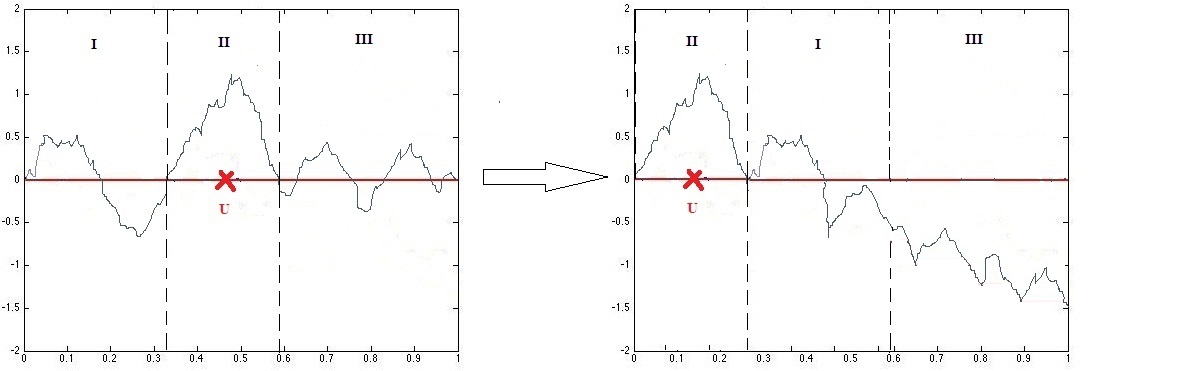}\\
Fig 4. Construction of Vervaat bridges via standard Brownian bridge.
\end{center}
\begin{proof} By Theorem \ref{PT}, $V(B^{\lambda,br})$ is entirely characterized by the distribution of the triple $(Z^{\lambda},B^{ex,Z^{\lambda}},F^{\lambda,1-Z^{\lambda}})$. It suffices to show that the distribution of the process defined as in \eqref{5} is determined by the same triple. Following from Perman et al. \cite[Theorem $3.1$]{PPY}, and Aldous and Pitman \cite[Corollary $5$]{AP}, conditioned on $\Delta:=D_U-G_U$, $(\widetilde{X}_t; 0 \leq t \leq \Delta)$ and $(\widetilde{X}_t; \Delta \leq t \leq 1)$ are independent, and $\Delta \stackrel{(d)}{=} Z^{\lambda}$ with density given by \eqref{aldous}.
Since $L_t^0(\widetilde{X})=0$ on $(0,\Delta)$, $(|\widetilde{X}(t)|-L_t^0(\widetilde{X}); 0 \leq t \leq \Delta)$ conditioned on $L_1^0(\widetilde{X})=|\lambda|$ is Brownian excursion of length $\Delta$, conditionally independent of $(|\widetilde{X}(t)|-L_t^0(\widetilde{X});  \Delta \leq t \leq 1)$ conditioned on $L_1^0(\widetilde{X})=|\lambda|$, which is a first passage bridge  of length $1-\Delta$ through $\lambda$ by the construction \eqref{4}.  
\end{proof}
\subsection{Vervaat bridges are not Markov}
\label{s32}
We ask if Vervaat bridges with non-zero endpoint are Markov processes. In the case of negative endpoint, the main problem is the lack of information on $Z^{\lambda}$: for $s \leq t \leq 1$, $V(B^{\lambda,br})_t$ depends not only on $V(B^{\lambda,br})_s$ but also on the event $\{Z^{\lambda} \leq s\}$. The following result gives a negative answer.
\begin{proposition}
\label{notmarkov}
$(V(B^{\lambda,br})_t; 0 \leq t \leq 1)$ where $\lambda<0$ is not Markov with respect to the induced filtration.
\end{proposition}

Before proving the proposition, we introduce some notations that are used in the current section, and in the rest of the paper.
For $x,y >0$, denote
$$
\tilde{q}_{t}(x,y):=\frac{1}{xy\sqrt{2 \pi t}} \left(\exp \left(-\frac{(x-y)^2}{2t}\right)- \exp \left(-\frac{(x+y)^2}{2t}\right)\right),
$$
$$
\tilde{q}_{t}(0,y)=\lim_{x\rightarrow 0^{+}}\tilde{q}_{t}(x,y)=
\dfrac{2}{\sqrt{2\pi t^{3}}}\exp\left(-\dfrac{y^{2}}{2t}\right)=\dfrac{2}{y}g_t(y) \quad \mbox{and} \quad \tilde{q}_{t}(0,0)=
\dfrac{2}{\sqrt{2\pi t^{3}}},
$$
where $g_t(y)$ is the density of the first hit to $y$ by Brownian motion given as in \eqref{Eq4}. Note that $\tilde{q}_{t}(x,y) y^{2}\,dy$ is the transition kernel of the three dimensional Bessel process. 
In the sequel, we aisely make use of conditioning and splicing Markov bridges, of which we refer to Fitzsimmons et al. \cite{FPY} for justifications.
\begin{proof}[Proof of Proposition \ref{notmarkov}] Fix $t_0 \in (0,1)$ and $x_0>0$. Let $T_{t_0}$ be the time of the first return by $V(B^{\lambda,br})$ to $0$ after time $t_0$. Consider the distribution of $T_{t_0}$ given $V(B^{\lambda,br})_{\frac{t_0}{2}}=0$ and $V(B^{\lambda,br})_{t_0}=x_0$. Given $T_{t_0}$, $(V(B^{\lambda,br})_t; t_0 \leq t \leq T_{t_0})$ and $(V(B^{\lambda,br})_t; T_{t_0} \leq t \leq 1)$ are two independent first passage bridges from $x_0>0$ to $0$, and respectively from $0$ to $\lambda<0$. Therefore, the conditional density of $T_{t_0}$ is given by
\begin{align}
\label{new1}
f_1(t)&=\frac{g_{t-t_0}(x_0) g_{1-t}(|\lambda|)}{g_{1-t_0}(x_0+|\lambda|)}1_{t > t_0} \notag\\
         &=\frac{C_1(t_0,x_0,\lambda)}{\sqrt{(t-t_0)^3(1-t)^3}} \exp \left( -\frac{x_0^2}{2(t-t_0)}-\frac{\lambda^2}{2(1-t)} \right) 1_{t>t_0},
\end{align}
for some $C_1(t_0,x_0,\lambda)>0$. Next we consider the distribution of $T_{t_0}$ given that $\forall u \in (0,t_0), V(B^{\lambda,br})_u>0$ and $V(B^{\lambda,br})_{t_0}=x_0$. Given $T_{t_0}$, $(V(B^{\lambda,br})_t; 0 \leq t \leq T_{t_0})$ is Brownian excursion, and the conditional density is computed via Bayes recipe:
\begin{align}
\label{new2}
f_2(t)&=C_2(t_0,x_0,\lambda) \frac{\tilde{q}_{t_0}(0,x_0) \tilde{q}_{t-t_0}(x_0,0)}{\tilde{q}_t(0,0)}x_0^2 f_{Z^{\lambda}}(t)1_{t>t_0} \notag\\
         &=C_2^{'}(t_0,x_0,\lambda) \frac{t}{\sqrt{(t-t_0)^3(1-t)^3}}\exp \left( -\frac{x_0^2}{2(t-t_0)}-\frac{\lambda^2}{2(1-t)} \right) 1_{t>t_0},
\end{align}
for some $C_2(t_0,x_0,\lambda)>0$ and $C_2^{'}(t_0,x_0,\lambda)>0$. 
Comparing \eqref{new1} to \eqref{new2}, we see that $f_2(t)=C_{1,2}(t_0,x_0,\lambda)tf_1(t)$ for some $C_{1,2}(t_0,x_0,\lambda)>0$. The two conditional densities of $T_{t_0}$ fail to be equal, which yields the desired result.  
\end{proof}

\quad It is well-known that the time reversal of any Markov process is still Markov. The result leads to the following corollary: 
\begin{corollary}
$(V(B^{\lambda,br})_t; 0 \leq t \leq 1))$ where $\lambda>0$ is not Markov with respect to the induced filtration.
\end{corollary}
\subsection{Semi-martingale decomposition of Vervaat bridges}
\label{s33}
The section is devoted to the canonical decomposition of Vervaat bridges with non-zero endpoint. Note that the splitting position in the case of negative endpoint (i.e. the first return to $0$) is a stopping time, while that of Vervaat bridges with positive endpoint (i.e. the last exit from $\lambda>0$) is a cooptional time but not a stopping time. The treatments in two cases are different. 
\subsubsection{Semi-martingale decomposition of $V(B^{\lambda,br})$ where $\lambda<0$}
\quad ~~Recall from Theorem \ref{PT} that $V(B^{\lambda,br})$ is split into Brownian excursion, followed by first passage bridge from $0$ to $\lambda$. The density of $Z^{\lambda}$ is given by \eqref{aldous}. We start by studying the canonical decomposition of $(V(B^{\lambda,br})_t;Z^{\lambda} \leq t \leq 1)$.
\begin{proposition}
\label{biay}
Let $\lambda<0$ and $V^{\lambda}:=V(B^{\lambda,br})$. Given $Z^{\lambda}:=\inf\{t >0; V(B^{\lambda,br})_t=0\}$,
$$\left( V^{\lambda}_t-\int_{Z^{\lambda}}^{t} \frac{1}{V^{\lambda}_s+\vert\lambda\vert}-\frac{ V^{\lambda}_s+\vert\lambda\vert}{1-Z^{\lambda}-s}ds \right)_{Z^{\lambda} \leq t \leq 1}$$
is Brownian motion.
\end{proposition}
\begin{proof} Conditional on the value of $Z^{\lambda}$, 
$(V^{\lambda}_t+\vert\lambda\vert)_{Z^{\lambda} \leq t \leq 1}$ is 
the three dimensional Bessel bridge from $|\lambda|$ to $0$.
This follows from Theorem \ref{PT}, and the following identity due to Biane and Yor \cite{BY}:
$$F^{\lambda,l}\stackrel{(d)}{=}BES(3)^{|\lambda| \rightarrow 0,l}+\lambda,$$ where $F^{\lambda,l}$ is first passage bridge of length $l$ from $0$ to $\lambda$, and $BES(3)^{|\lambda| \rightarrow 0,l}$ is the three dimensional Bessel bridge from $|\lambda|$ to $0$. 
The proposition follows from the semi-martingale decomposition of Bessel bridges, see e.g. Revuz and Yor \cite[Chapter XI.$3$)]{RY}.
\end{proof}

We deal with the canonical decomposition of $ (V(B^{\lambda,br})_{t \wedge Z^{\lambda}};0 \leq t \leq 1)$. The process is Brownian excursion of length $Z^{\lambda}$, absorbed at $0$ after time $Z^{\lambda}$.
\begin{proposition}
\label{keylem}
Let $\lambda<0$ and
$$J_{t}^{\lambda}(y):=\int_{t}^1 \frac{\tilde{q}_{s-t}(0,y)}{\tilde{q}_s(0,0)}f_{Z^{\lambda}}(s) ds,$$
$$\mathring{J}_{t}^{\lambda}(y):=\int_{t}^1 \frac{1}{\sqrt{s-t}} \frac{\tilde{q}_{s-t}(0,y)}{\tilde{q}_s(0,0)}f_{Z^{\lambda}}(s) ds.$$
Let $V^{\lambda}:=V(B^{\lambda,br})$, then
\begin{displaymath}
(Y^{\lambda}_{t})_{0 \leq t \leq 1}:=\left(V^{\lambda}_{t \wedge Z^{\lambda}}-\int_{0}^{t\wedge Z^{\lambda}}\dfrac{ds}{V^{\lambda}_{s}}+
\int_{0}^{t\wedge Z^{\lambda}}
\dfrac{V^{\lambda}_{s}\mathring{J}^{\lambda}_{s}(V^{\lambda}_{s})}{J^{\lambda}_{s}(V^{\lambda}_{s})}
\,ds\right)_{0 \leq t \leq 1}
\end{displaymath}
is Brownian motion with respect to the filtration of $V(B^{\lambda,br})$, stopped at time $Z^{\lambda}$.
\end{proposition}
\begin{proof} Let $\varepsilon\in(0,1)$. Introduce $(B^{\lambda,\varepsilon}_{t};t\geq 0)$ Brownian motion with the starting point $B^{\lambda,\varepsilon}_{0} \stackrel{(d)}{=}V(B^{\lambda,br})_{\varepsilon \wedge Z^{\lambda}}$. The density of this distribution is:
\begin{displaymath}
\mu_{\varepsilon}^{\lambda}(x)=\tilde{q}_{\varepsilon}(0,x) x^2 J_{\varepsilon}^{\lambda}(x) \quad \mbox{for}~x>0 \quad \mbox{and}\quad \mu_{\varepsilon}^{\lambda}(0)=\int_0^{\varepsilon} f_{Z^{\lambda}}(s)ds.
\end{displaymath}
Let $T^{\lambda,\varepsilon}_{0}$ be the first time at which $B^{\lambda,\varepsilon}$ hits $0$. For $\varepsilon \leq t \leq 1$, the distribution of $(V(B^{\lambda,br})_{s \wedge Z^{\lambda}};\varepsilon \leq s \leq t)$ is absolutely continuous with respect to that of $(B^{\lambda,\varepsilon}_{(s-\varepsilon)\wedge T^{\lambda,\varepsilon}_{0}};\varepsilon\leq s\leq t)$. In fact, conditionally on $Z^{\lambda}>t$ and the value of $V(B^{\lambda,br})_t$, $(V(B^{\lambda,br})_{s \wedge Z^{\lambda}};\varepsilon \leq s \leq t)$ has the same distribution as a three dimensional Bessel bridge. It is the same for $(B^{\lambda,\varepsilon}_{(s-\varepsilon)\wedge T^{\lambda,\varepsilon}_{0}};\varepsilon\leq s\leq t)$ conditionally on $ T^{\lambda,\varepsilon}_{0}>t-\varepsilon$ and the value of $B^{\lambda,\varepsilon}_{t-\varepsilon}$. The corresponding density 
$D_t^{\lambda,\varepsilon}$ equals
\begin{align}
&~\quad 1_{B_0^{\lambda,\varepsilon}=0}+1_{B_0^{\lambda,\varepsilon}>0,T_0^{\lambda,\varepsilon} \leq t-\varepsilon}+1_{T_0^{\lambda,\varepsilon} > t-\varepsilon} \frac{\tilde{q}_{\varepsilon}(0,B^{\lambda,\varepsilon}_{0}) (B^{\lambda,\varepsilon}_{0})^2
\tilde{q}_{t-\varepsilon}(B^{\lambda,\varepsilon}_{0},B^{\lambda,\varepsilon}_{t-\varepsilon}) (B^{\lambda,\varepsilon}_{t-\varepsilon})^2 J_t^{\lambda}(B_{t-\varepsilon}^{\lambda,\varepsilon})}{\mu_{\varepsilon}^{\lambda}(B^{\lambda,\varepsilon}_{0})B^{\lambda,\varepsilon}_{0}B^{\lambda,\varepsilon}_{t-\varepsilon}\tilde{q}_{t-\varepsilon}(B^{\lambda,\varepsilon}_{0},B^{\lambda,\varepsilon}_{t-\varepsilon})} \notag\\
                                            &=  1_{B_0^{\lambda,\varepsilon}=0}+1_{B_0^{\lambda,\varepsilon}>0,T_0^{\lambda,\varepsilon} \leq t-\varepsilon}+1_{T_0^{\lambda,\varepsilon} > t-\varepsilon}\frac{B^{\lambda,\varepsilon}_{t-\varepsilon} J^{\lambda}_t(B^{\lambda,\varepsilon}_{t-\varepsilon})}{B^{\lambda,\varepsilon}_0 J^{\lambda}_{\varepsilon}(B^{\lambda,\varepsilon}_0)}. \notag
\end{align}
In addition, $\dfrac{\partial J_{t}^{\lambda}(y)}{\partial y}=-y\mathring{J}_{t}^{\lambda}(y)$ for $t\in(0,1)$. We get the quadratic covariation
\begin{displaymath}
d\left[\log(D_{t\wedge T_0^{\lambda,\varepsilon}}^{\lambda,\varepsilon}),
B_{(t-\varepsilon)\wedge T_0^{\lambda,\varepsilon}}
^{\lambda,\varepsilon}\right]=1_{T_0^{\lambda,\varepsilon} > t-\varepsilon}
\dfrac{J^{\lambda}_t(B^{\lambda,\varepsilon}_{t-\varepsilon})-
(B^{\lambda,\varepsilon}_{t-\varepsilon})^{2} \mathring{J}^{\lambda}_t(B^{\lambda,\varepsilon}_{t-\varepsilon})}
{B^{\lambda,\varepsilon}_{t-\varepsilon} J^{\lambda}_t(B^{\lambda,\varepsilon}_{t-\varepsilon})}dt,
\end{displaymath}
where $[\cdot,\cdot]$ is the quadratic covariation of two processes. Applying Girsanov's theorem, we obtain that $(Y^{\lambda}_{t})_{t\geq \varepsilon}$ is a continuous martingale relative to the filtration of $(V(B^{\lambda,br})_t; \varepsilon \leq t \leq 1)$ with quadratic variation $(t-\varepsilon)\wedge
(Z^{\lambda}-\varepsilon)^{+}$. Since this holds for all $\varepsilon$ sufficiently small, we have proved the proposition.  
\end{proof}

Following from Proposition \ref{biay} and Proposition \ref{keylem}, we obtain:
\begin{theorem}
Let $\lambda<0$ and $V^{\lambda}:=V(B^{\lambda,br})$. Then
$$
\Bigg(V^{\lambda}_{t}-\int_{0}^{t\wedge Z^{\lambda}}\dfrac{ds}{V^{\lambda}_{s}}+
\int_{0}^{t\wedge Z^{\lambda}}\dfrac{V^{\lambda}_{s}\mathring{J}_{s}^{\lambda}(V^{\lambda}_{s})}{J_{s}^{\lambda}(V^{\lambda}_{s})}
\,ds -\int_{Z^{\lambda}}^{t} \frac{1}{ V^{\lambda}_s+\vert\lambda\vert}-
\frac{ V^{\lambda}_s+\vert\lambda\vert}{1-Z^{\lambda}-s}ds\Bigg)_{0\leq t\leq 1}
$$
is Brownian motion.
\end{theorem}
\subsubsection{Semi-martingale decomposition of $V(B^{\lambda,br})$ where $\lambda>0$}
\quad~~Recall from Corollary \ref{PT2} that $V(B^{\lambda,br})$ is decomposed into the three dimensional Bessel bridge from $0$ to $\lambda$, followed by Brownian excursion above $\lambda$. The density of $\widehat{Z}^{\lambda}$ is given by $f_{\widehat{Z}^{\lambda}}(t)=f_{Z^{\lambda}}(1-t)$, as in \eqref{aldous}.

For $x,y\geq 0$, let $\mathbb{Q}^{t}_{x,y}$ be the distribution of the three dimensional Bessel bridge of length $t$ from $x$ to $y$. Let $(R_{t})_{t\geq 0}$ be the three dimensional Bessel process starting from $0$. The key idea is to show that for any $t\in[0,1)$, the distribution of $(V(B^{\lambda,br})_{s}; 0 \leq s \leq t)$ is absolutely continuous with respect to that of $(R_{s}; 0 \leq s \leq t)$, identify the corresponding density $D^{\lambda}_{t}$, and deduce by applying Girsanov's theorem the canonical decomposition of $V(B^{\lambda,br})$.

Let
$\theta^{\lambda}_{t}:=\sup\lbrace s\in[0,t]; R_{s}\leq \lambda\rbrace$.
Note that if $R_{t}\leq\lambda$ then $\theta^{\lambda}_{t}=t$. We start with a lemma computing the joint distribution of $(R_t,\theta_t)$ in the case where the last exit from $\lambda$ has not been occurred.

\begin{lemma}
\label{Lem1}
On the event $R_{t}>\lambda$, the joint distribution of $(R_{t},\theta^{\lambda}_{t})$ is:
$$
\tilde{q}_{t}(0,y)\dfrac{g_{t-s}(y-\lambda)g_s(\lambda)}{g_t(y)}\,1_{0<s<t}\,ds\,
1_{y>\lambda}y^{2}\,dy.
$$
Conditionally on $R_{t}>\lambda$, the value of $R_{t}$ and of $\theta^{\lambda}_{t}$, $(R_{s};0\leq s\leq \theta^{\lambda}_{t})$ and $(R_{\theta^{\lambda}_{t}+s}-\lambda;0\leq s\leq t- \theta^{\lambda}_{t})$ are independent and follow the law $\mathbb{Q}^{\theta^{\lambda}_{t}}_{0,\lambda}$ respectively $\mathbb{Q}^{t-\theta^{\lambda}_{t}}_{0,R_{t}-\lambda}$.
\end{lemma}
\begin{proof}
Let $y>\lambda$. Conditionally on $R_{t}=y$, $(R_{t-s};0\leq s\leq t)$ is first passage bridge from $y$ to $0$, and $t-\theta^{\lambda}_{t}$ is the first time where it hits $\lambda$. Thus conditionally on $R_{t}=y$, $t-\theta^{\lambda}_{t}$ is distributed as
\begin{displaymath}
\dfrac{g_s(y-\lambda)g_{t-s}(\lambda)}{g_t(y)}\,1_{0<s<t}\,ds.
\end{displaymath}
Moreover, conditionally on $R_{t}=y$ and the value of $\theta^{\lambda}_{t}$, $(R_{t-s};0\leq s\leq t-\theta^{\lambda}_{t})$ and $(R_{\theta^{\lambda}_{t}-s};0\leq s\leq \theta^{\lambda}_{t})$ are two independent first passage bridges, from $y$ to $\lambda$ and from $\lambda$ to $0$.  
\end{proof}

In the next proposition, we establish the absolute continuity between the distribution of $(V(B^{\lambda,br})_{s};0\leq s\leq t)$ and that of $(R_{s};0\leq s\leq t)$. The Radon-Nikodym density $D^{\lambda}_{t}$ is expressed as a deterministic function of $t, R_{t}$ and $\theta^{\lambda}_{t}$.
\begin{proposition}
For any $t\in[0,1)$, the distribution of $(V(B^{\lambda,br})_{s};0\leq s\leq t)$ is absolutely continuous with respect to that of $(R_{s}; 0 \leq s \leq t)$. The corresponding density is:
\begin{align}
D^{\lambda}_{t}&=\int_{t}^{1}\dfrac{\tilde{q}_{s-t}(R_{t},\lambda)}{\tilde{q}_{s}(0,\lambda)}f_{\widehat{Z}^{\lambda}}(s)\,ds+
1_{R_{t}>\lambda}\dfrac{(1-\theta^{\lambda}_{t})(R_{t}-\lambda)}{\sqrt{(1-t)^{3}}R_{t}}\exp\left(\dfrac{\lambda^{2}}{2}\right)
\exp\left(-\dfrac{(R_{t}-\lambda)^{2}}{2(1-t)}\right) \notag\\
   &:=\Phi^{\lambda}(t,R_t,\theta^{\lambda}_t). \label{Bigphi}
\end{align}
In particular, $(D^{\lambda}_{t};0\leq t<1)$ is continuous and there is no discontinuity as $R_{t}$ crosses $\lambda$. 
\end{proposition}
\begin{proof} Write
$D^{\lambda}_{t}=D^{1,\lambda}_{t}+D^{2,\lambda}_{t}$, where
$D^{1,\lambda}_{t}:=1_{\widehat{Z}^{\lambda}>t}D^{\lambda}_{t}$ and $D^{2,\lambda}_{t}:=1_{\widehat{Z}^{\lambda} <t} D^{\lambda}_{t}$. On the event $R_{t}<\lambda$, we have $D^{\lambda}_{t}=D^{1,\lambda}_{t}$. Conditionally on $\widehat{Z}^{\lambda}>t$ and the value of $V(B^{\lambda,br})_{t}$, $(V(B^{\lambda,br})_{s};0\leq s\leq t)$ is the three dimensional Bessel bridge from $0$ to $V(B^{\lambda,br})_{t}$, i.e. the same distribution as that of $(R_{s};0\leq s\leq t)$ conditioned on the value of $R_{t}$. Conditionally on $\widehat{Z}^{\lambda}>t$ and the position of $\widehat{Z}^{\lambda}$, $V(B^{\lambda,br})_{t}$ is distributed as
\begin{displaymath}
\dfrac{\tilde{q}_{t}(0,y)\tilde{q}_{\widehat{Z}^{\lambda}-t}(y,\lambda)}{\tilde{q}_{\widehat{Z}^{\lambda}}(0,\lambda)}\,1_{y>0}\,y^{2}dy.
\end{displaymath}
Therefore,
$$
D^{1,\lambda}_{t}=\int_{t}^{1}\dfrac{\tilde{q}_{s-t}(R_{t},\lambda)}{\tilde{q}_{s}(0,\lambda)}f_{\widehat{Z}^{\lambda}}(s)\,ds .$$
Conditionally on $\widehat{Z}^{\lambda}<t$, the position of $\widehat{Z}^{\lambda}$ and the value of $V(B^{\lambda,br})_{t}$, $(V(B^{\lambda,br})_{s};0\leq s\leq \widehat{Z}^{\lambda})$ and $(V(B^{\lambda,br})_{\widehat{Z}^{\lambda}+s}-\lambda;0\leq s\leq t-\widehat{Z}^{\lambda})$ are independent and follow the distribution $\mathbb{Q}^{\widehat{Z}^{\lambda}}_{0,\lambda}$ respectively $\mathbb{Q}^{t-\widehat{Z}^{\lambda}}_{0,V(B^{\lambda,br})_{t}-\lambda}$. These are the same conditional distributions as in Lemma \ref{Lem1}. On the event $\widehat{Z}^{\lambda}<t$, the joint distribution of $(V(B^{\lambda,br})_{t},\widehat{Z}^{\lambda})$ is given by
\begin{displaymath}
f_{\widehat{Z}^{\lambda}}(s)\dfrac{\tilde{q}_{t-s}(0,y-\lambda)\tilde{q}_{1-t}(y-\lambda,0)}{\tilde{q}_{1-s}(0,0)}\,1_{y>\lambda}\,(y-\lambda)^{2}dy\,1_{0<s<t}\,ds.
\end{displaymath}
We have then,
\begin{equation*}
\begin{split}
D^{2,\lambda}_{t}=&1_{R_{t}>\lambda}
\dfrac{f_{\widehat{Z}^{\lambda}}(\theta^{\lambda}_{t})\dfrac{\tilde{q}_{t-\theta^{\lambda}_{t}}(0,R_{t}-\lambda)\tilde{q}_{1-t}(R_{t}-\lambda,0)}{\tilde{q}_{1-\theta^{\lambda}_{t}}(0,0)}(R_{t}-\lambda)^{2}}
{\tilde{q}_{t}(0,R_{t})\dfrac{g_{t-\theta^{\lambda}_t}(R_t-\lambda)g_{\theta^{\lambda}_t}(\lambda)}{g_t(R_t)}\,R_{t}^{2}}\\
=&1_{R_{t}>\lambda}\dfrac{(1-\theta^{\lambda}_{t})(R_{t}-\lambda)}{\sqrt{(1-t)^{3}}R_{t}}\exp\left(\dfrac{\lambda^{2}}{2}\right)
\exp\left(-\dfrac{(R_{t}-\lambda)^{2}}{2(1-t)}\right).
\end{split}
\end{equation*}
\end{proof}

\begin{lemma}
\label{LemIntegrale}
For any $t\in(0,1)$ and $a\geq 0$:
\begin{displaymath}
\int_{t}^{1}\dfrac{ds}{\sqrt{(1-s)(s-t)}}\exp\left(-\dfrac{a}{s-t}\right)=\sqrt{\pi}
\int^{+\infty}_{\frac{a}{1-t}}e^{-u}\dfrac{du}{\sqrt{u}}.
\end{displaymath}
\end{lemma}
\begin{proof}
By change of variables $z:=\dfrac{1-s}{1-t}$, we obtain:
\begin{displaymath}
\int_{t}^{1}\dfrac{ds}{\sqrt{(1-s)(s-t)}}\exp\left(-\dfrac{a}{s-t}\right)=
\int_{0}^{1}\dfrac{dz}{\sqrt{z(1-z)}}\exp\left(-\dfrac{a}{(1-t)z}\right):=
\varphi(\frac{a}{1-t}),
\end{displaymath}
Note that
\begin{displaymath}
\varphi(x):=\int_{0}^{1}\dfrac{dz}{\sqrt{z(1-z)}}\exp\left(-\dfrac{x}{z}\right)\stackrel{(\#)}{=}\int_{1}^{+\infty}\dfrac{dv}{v\sqrt{v-1}}e^{-xv},
\end{displaymath}
where $(\#)$ is obtained by taking $z=v^{-1}$. Following from Gradshteyn and Ryzhik \cite[$3.363 (2)$]{GR}, we have $
\varphi(x)=\sqrt{\pi}\int_{x}^{+\infty}e^{-u}\dfrac{du}{\sqrt{u}}$ .  
\end{proof}

Let
\begin{displaymath}
\Phi^{1,\lambda}(t,y):=\dfrac{1}{2\sqrt{2}y}
\exp\left(\frac{\lambda^{2}}{2}\right)\int_{\frac{(y-\lambda)^{2}}{2(1-t)}}^{\frac{(y+\lambda)^{2}}{2(1-t)}} 
e^{-u}\dfrac{du}{\sqrt{u}},
\end{displaymath}
\begin{displaymath}
\Phi^{2,\lambda}(t,y):=\dfrac{(y-\lambda)}{\sqrt{(1-t)^{3}}y}\exp\left(\dfrac{\lambda^{2}}{2}\right)
\exp\left(-\dfrac{(y-\lambda)^{2}}{2(1-t)}\right).
\end{displaymath}
According to the formula \eqref{Bigphi} and Lemma \ref{LemIntegrale}:
\begin{equation}
\label{phi}
\Phi^{\lambda}(t,y,\theta)=\Phi^{1,\lambda}(t,y)+(1-\theta)(0\vee \Phi^{2,\lambda}(t,y)).
\end{equation}
Observe that $\Phi^{2,\lambda}$ is of class $\mathcal{C}^{1}$. $\Phi^{1,\lambda}$ and the partial derivative $\partial_{1}\Phi^{1,\lambda}$ are continuous as functions of $(t,y)$. However, $\partial_{2}\Phi^{1,\lambda}(t,y)$ is not defined at $y=\lambda$:
\begin{displaymath}
\partial_{2}\Phi^{1,\lambda}(t, \lambda^{+})-\partial_{2}\Phi^{1,\lambda}(t, \lambda^{-})=
-\dfrac{1}{\sqrt{1-t}\lambda}\exp\left(\dfrac{\lambda^{2}}{2}\right),
\end{displaymath}
\begin{align}
\partial_{2}\Phi^{\lambda}(t, \lambda^{+},\theta)-\partial_{2}\Phi^{\lambda}(t,\lambda^{-},\theta)
&=\partial_{2}\Phi^{1,\lambda}(t, \lambda^{+})-\partial_{2}\Phi^{1,\lambda}(t, \lambda^{-})
+(1-\theta)\partial_{2}\Phi^{2,\lambda}(t,\lambda) \notag\\
&=\dfrac{(t-\theta)}{\sqrt{(1-t)^{3}}\lambda}\exp\left(\dfrac{\lambda^{2}}{2}\right). \notag
\end{align}

For $t>0$, let
$
W_{t}:=R_{t}-\int_{0}^{t}\dfrac{ds}{R_{s}},
$
where $(W_{t};t\geq 0)$ is Brownian motion starting from $0$, predictable with respect to the filtration of $(R_{t};t\geq 0)$.

\begin{lemma}
\label{LemMartDec}
For all $t\in[0,1)$ and $\lambda>0$,
\begin{displaymath}
D^{\lambda}_{t}=1+\int_{0}^{t}\partial_{2}\Phi^{\lambda}(s,R_{s},\theta^{\lambda}_{s})\,dW_{s}.
\end{displaymath}
\end{lemma}
\begin{proof}
Remark that we cannot apply directly It\^{o}'s formula to $\Phi^{\lambda}(t,R_{t},\theta^{\lambda}_{t})$, since $\Phi^{\lambda}$ is not regular enough.
It is easy to check that $\Phi^{2,\lambda}$ and $\Phi^{1,\lambda}$ outside $\lbrace y=\lambda\rbrace$ satisfy the PDE:
\begin{displaymath}
\dfrac{1}{2}\partial_{2,2}\Phi(t,y)+\dfrac{1}{y}\partial_{2}\Phi(t,y)+\partial_{1}\Phi(t,y)=0.
\end{displaymath}
Let $(L^{\lambda}_{t}(R);t \geq 0)$ be the local times at level $\lambda$ of $(R_{t};t\geq 0)$. Applying It\^{o}-Tanaka's formula, and taking into account the discontinuity of the partial derivative $\partial_{2}$ at level $y=\lambda$, we get:
\begin{displaymath}
\Phi^{1,\lambda}(t,R_{t})=1+\int_{0}^{t}\partial_{2}\Phi^{1,\lambda}(s,R_{s})\,dW_{s}-
\dfrac{1}{\lambda}\exp\left(\dfrac{\lambda^{2}}{2}\right)\int_{0}^{t}\dfrac{1}{\sqrt{1-s}}
\,dL^{\lambda}_{s}(R).
\end{displaymath}
\begin{displaymath}
0\vee \Phi^{2,\lambda}(t,R_{t})=\int_{0}^{t}1_{R_{s}>\lambda}
\partial_{2}\Phi^{2,\lambda}(s,R_{s})\,dW_{s}+
\dfrac{1}{\lambda}\exp\left(\dfrac{\lambda^{2}}{2}\right)\int_{0}^{t}\dfrac{1}{\sqrt{(1-s)^{3}}}
\,dL^{\lambda}_{s}(R).
\end{displaymath}
The process $(1-\theta^{\lambda}_{t})$ is not continuous but it is constant on the intervals of time where $0\vee \Phi^{2,\lambda}(t,R_{t})$ is positive. According to the derivation rule in Revuz and Yor \cite[Theorem $4.2$, Chapter VI]{RY},
\begin{align*}
(1-\theta^{\lambda}_{t})(0\vee \Phi^{2,\lambda}(t,R_{t}))
&=\int_{0}^{t}(1-\theta^{\lambda}_{s})
d(0\vee \Phi^{2,\lambda}(s,R_{s}))\\
&=\int_{0}^{t}1_{R_{s}>\lambda}
(1-\theta^{\lambda}_{s})\partial_{2}\Phi^{2,\lambda}(s,R_{s})\,dW_{s}+
\dfrac{1}{\lambda}e^{\frac{\lambda^2}{2}}\int_{0}^{t}\dfrac{(1-\theta^{\lambda}_{s})}
{\sqrt{(1-s)^{3}}}\,dL^{\lambda}_{s}(R)\\
&=\int_{0}^{t}1_{R_{s}>\lambda}
(1-\theta^{\lambda}_{s})\partial_{2}\Phi^{2,\lambda}(s,R_{s})\,dW_{s}+
\dfrac{1}{\lambda}e^{\frac{\lambda^2}{2}}\int_{0}^{t}\dfrac{1}
{\sqrt{1-s}}\,dL^{\lambda}_{s}(R).
\end{align*}
on the support of $dL^{\lambda}_{s}(R)$, $(1-\theta^{\lambda}_{s})$ being equal to $1-s$. Finally
\begin{multline*}
\Phi^{1,\lambda}(t,R_{t})+(1-\theta^{\lambda}_{t})(0\vee \Phi^{2,\lambda}(t,R_{t}))\\=
1+\int_{0}^{t}(\partial_{2}\Phi^{1,\lambda}(s,R_{s})+(1-\theta^{\lambda}_{s})\partial_{2}
\Phi^{2,\lambda}(s,R_{s}))1_{R_s > \lambda}\,dW_{s},
\end{multline*}
which leads to the desired result.  
\end{proof}

\begin{theorem}
Let $\lambda>0$ and for $t\in(0,1)$,
\begin{displaymath}
\tilde{\theta}_{t}^{\lambda,br}:=\sup\lbrace s\in[0,t]\vert V(B^{\lambda, br})_{s}\leq \lambda\rbrace.
\end{displaymath}
Let $V^{\lambda}:=V(B^{\lambda,br})$, then
\begin{displaymath}
\left(V^{\lambda}_{t}-\int_{0}^{t}\dfrac{ds}{V^{\lambda}_{s}}-\int_{0}^{t}\dfrac{\partial_{2}\Phi^{\lambda}}{\Phi^{\lambda}}
(s,V^{\lambda}_{s},\tilde{\theta}_{s}^{\lambda,br})\,ds\right)_{0\leq t\leq 1}
\end{displaymath}
is Brownian motion.
\end{theorem}
\begin{proof}
For $t\in[0,1)$, let
$$
X_{t}:=V(B^{\lambda,br})_{t}-\int_{0}^{t}\dfrac{ds}{V(B^{\lambda,br})_{s}}.
$$
The distribution of $(X_{s};0\leq s\leq t)$ is absolutely continuous with respect to that of $(W_{s};0\leq s\leq t)$, with Radon-Nikodym density $D^{\lambda}_{t}$. By Lemma \ref{LemMartDec}, we get the quadratic covariation
\begin{displaymath}
[\log(D^{\lambda}),W]_{t}=
\int_{0}^{t}\dfrac{\partial_{2}\Phi^{\lambda}}{\Phi^{\lambda}}(s,R_{s},\theta^{\lambda}_{s})\,ds.
\end{displaymath}
Applying Girsanov's theorem, we obtain that
\begin{displaymath}
X_{t}-\int_{0}^{t}\dfrac{\partial_{2}\Phi^{\lambda}}{\Phi^{\lambda}}
(s,V(B^{\lambda,br})_{s},\tilde{\theta}_{s}^{\lambda,br})\,ds
\end{displaymath}
is Brownian motion.  
\end{proof}
\subsection{Relation with drifting excursion}
\label{s34}
Bertoin \cite{Bertoin} studied a fragmentation process by considering normalized Brownian excursion dragged down by drift $\lambda<0$:
$$B_t^{ex, \lambda \downarrow}:=B^{ex}_t+\lambda t, ~\mbox{for}~0 \leq t \leq 1.$$

Note that $V(B^{\lambda,br})$ where $\lambda<0$ looks similar to this process, except that $B^{ex, \lambda \downarrow}$ always stays above the line $t \rightarrow \lambda t$, while Vervaat bridges do not share this property. An interesting question is whether conditioned on staying above the dragging line, the Vervaat bridge is absolutely continuous with respect to drifting excursion. 
To this end, we need to justify that the conditioning event has positive probability, as shown in the next proposition.
\begin{proposition}
\label{comp}
Let $\lambda<0$. Then
$$ \mathbb{P}(\forall t \in (0,1),~V(B^{\lambda,br})_t >\lambda t)=1-|\lambda| \exp\left(\frac{\lambda^2}{2}\right) \int_{|\lambda|}^{\infty} \exp\left(-\frac{t^2}{2}\right)dt.$$
\end{proposition}
\begin{proof} According to Schweinsberg \cite[Proposition $15$]{Schweinsberg}, fix $x \in [\lambda,0]$, the probability for first passage bridge through $\lambda$ to stay above the dragging line tying $x$ to $\lambda$ is:
\begin{equation}
\label{15}
\mathbb{P}(\forall t \in [0,l],~F^{\lambda,l}(t)>x-(x-\lambda)t)=\frac{|x|}{|\lambda|}.
\end{equation}
Therefore,
\begin{align}
\mathbb{P}(\forall t \in (0,1),~V(B^{\lambda,br})_t >\lambda t) &=\int_0^1  \mathbb{P}\left(\forall s \in (t,1),~V(B^{\lambda,br})_s >\lambda s | Z^{\lambda}=t\right)f_{Z^{\lambda}}(t)dt  \notag\\
                                                                                                    &\stackrel{(*)}{=}\int_0^1  t \frac{|\lambda|}{\sqrt{2 \pi t(1-t)^3}} \exp\left(-\frac{\lambda^2 t}{2(1-t)}\right)dt  \notag\\
                                                                                                    &=\mathbb{E}Z^{\lambda},   \notag 
\end{align}
where $(*)$ follows from \eqref{15}. It suffices to apply Pitman \cite[Lemma $4.10$]{Pitman}:
$$\mathbb{E}Z^{\lambda}=1-\lambda \exp\left(\frac{\lambda^2}{2}\right) \int_{\lambda}^{\infty} \exp\left(-\frac{t^2}{2}\right)dt.$$
\end{proof}

We know that the Vervaat bridge with negative endpoint conditioned to stay above the dragging line is well-defined. Moreover, the distribution of the first return to $0$ is given by
\begin{equation}
\label{16}
f_{\widetilde{Z}^{\lambda}}(t)=\frac{t}{1-|\lambda| \exp\left(\frac{\lambda^2}{2}\right) \int_{|\lambda|}^{\infty} \exp\left(-\frac{t^2}{2}\right)dt}f_{Z^{\lambda}}(t). 
\end{equation}
\begin{corollary}
Let $\lambda<0$. Given $\widetilde{Z}^{\lambda}$ the time of the first return to $0$ by $(V(B^{\lambda,br})_t; 0 \leq t \leq 1|\forall t \in (0,1), V(B^{\lambda,br})_t>\lambda t)$, whose distribution density is given by \eqref{16}, the path is decomposed into two (conditionally) independent pieces:
\begin{itemize}
\item
$\left(V(B^{\lambda,br})_u;0 \leq u \leq \widetilde{Z}^{\lambda}|\forall t \in (0,1), V(B^{\lambda,br})_t>\lambda t\right)$ is excursion of length $\widetilde{Z}^{\lambda}$;
\item
 $\left(V(B^{\lambda,br})_u;\widetilde{Z}^{\lambda} \leq u \leq 1|\forall t \in (0,1), V(B^{\lambda,br})_t>\lambda t\right)$ is first passage bridge of length $1-\widetilde{Z}^{\lambda}$ conditioned to stay above $t \rightarrow \lambda(t+\widetilde{Z}^{\lambda})$ for $t \in (0,1-\widetilde{Z}^{\lambda})$.
\end{itemize}
In addition, $(V(B^{\lambda,br})_t; 0 \leq t \leq 1|\forall t \in (0,1), V(B^{\lambda,br})_t>\lambda t)$ is absolutely continuous with respect to $(B_t^{ex,\lambda \downarrow}; 0\leq t \leq 1)$. The corresponding density is: $$\frac{H}{1-|\lambda| \exp\left(\frac{\lambda^2}{2}\right) \int_{|\lambda|}^{\infty} \exp\left(-\frac{t^2}{2}\right)dt},$$
where $H:=\inf \{t>0; B^{ex, \lambda \downarrow}_{t}<0\}$.
\end{corollary}
\begin{proof} Following from Bertoin \cite[Proposition $11$]{Bertoin}, $H$ is distributed as in \eqref{aldous}. According to Chassaing and Janson \cite[Theorem $2.6$]{CJ}, conditioned on $H$, $(B^{ex, \lambda \downarrow}_t; 0 \leq t \leq H)$ is Brownian excursion of length $H$. Moreover, Schweinsberg \cite[Proposition $4$]{Schweinsberg} states that given $H$, $(B^{ex, \lambda \downarrow}_t; H \leq t \leq 1)$ is a first passage bridge of length $1-H$ conditioned to stay above the line $t \rightarrow \lambda(t+H)$ for $t \in (0,1-H)$, (conditionally) independent of the excursion. By change of measures, we obtain the same triple characterization in distribution. 
\end{proof}
\subsection{Convex minorant of Vervaat bridges}
\label{s35}
In this subsection, we study some properties of the convex minorant of Vervaat bridges $V(B^{\lambda,br})$ where $\lambda<0$. The convex minorant of a real-valued function $(X_t; 0 \leq t \leq 1)$ is the maximal convex function $(C_t; t \in [0,1])$ such that $ \forall t \in [0,1], C_t \leq X_t $. We refer to the points where the convex minorant equals the process as vertices. Note that these points are also the endpoints of the linear segments. We refer readers to Pitman and Ross \cite{PR}, and Abramson et al. \cite{Abr} for background.

Similar to the computation in Proposition \ref{comp} , we have an explicit formula for the distribution of the last segment's slopes.
\begin{corollary}
Denote $s_l$ the slope of the last segment of the convex minorant of $V(B^{\lambda,br})$. For $a \in [\lambda,0]$, we have
$$\mathbb{P}(s_l \in [\lambda, a])=1+a \exp\left(\frac{\lambda^2}{2}\right) \int_{|\lambda|}^{\infty} \exp\left(-\frac{t^2}{2}\right)dt.$$
\end{corollary}

As discussed in Pitman and Ross \cite{PR}, first passage bridge of length $1$ can only have accumulations of linear segments at its start point. As mentioned before, the greatest difference between Vervaat bridges and first passage bridges of length $1$ is that the former starts with a piece of excursion, while the latter returns to $0$ immediately. We expect that almost surely, Vervaat bridges have a finite number of segments.
\begin{proposition}
The number of segments of the convex minorant of  $V(B^{\lambda,br})$ where $\lambda<0$ is almost surely finite.
\end{proposition}
\begin{proof} Consider a sample path of Brownian bridge $B^{\lambda,br}$ where $\lambda<0$ and $1-A^{\lambda}:=\argmin B^{\lambda,br}$, which is almost surely unique. Note that $V(B^{\lambda,br})_t >0$ for $t \in (0,A^{\lambda}]$. Consequently, the first vertex of the Vervaat bridge $\alpha_1>A^{\lambda}$ almost surely. By Pitman and Ross \cite{PR}, there can be only a finite number of segments on $[\alpha_1,1]$, since accumulations can only happen at $0$ in the path of $B^{\lambda,br}$ restricted to $[0,1-A^{\lambda}]$. Thus, the number of segments of the Vervaat bridge is almost surely finite.  
\end{proof}
\section{The Vervaat transform of Brownian motion}
\label{s4}
In this section, we study the Vervaat transform of Brownian motion. We first prove that the process is not a Markov process. Next, $V(B)$ is shown to be semi-martingale, and its canonical decomposition is given. The computation is essentially based on the results of Subsection \ref{s33}. Finally, we compute the mean and the variance of this process.
\subsection{$V(B)$ is not Markov}
\label{s41}
A crucial property of $V(B)$ is that $V(B)_1=B_1$. We encounter two cases. If $B_1>0$, then $V(B)$ never returns to $0$ and stays positive along the path. Otherwise $V(B)_1 =B_1 \leq 0$. By path continuity, $V(B)$ has to hit $0$ somewhere in the path.
\begin{proposition}
$(V(B)_t; 0 \leq t \leq 1)$ is not Markov.
\end{proposition}
\begin{proof} Fix $x_0>0$. Given $V(B)_{\frac{1}{4}}=0$ and $V(B)_{\frac{1}{2}}=x_0$, the conditional distribution of $V(B)_1$ is supported on the negative half line $(-\infty,0]$, since once $V(B)$ hits $0$ in the path, it has to end negatively. On the other hand, given $V(B)_t>0$ for all $t \in (0,\frac{1}{2}]$ and $V(B)_{\frac{1}{2}}=x_0$, the support of $V(B)_1$ is clearly the whole real line. 
These two conditional distributions fail to be equal, which yields the desired result.  
\end{proof} 

In other words, $\{\widetilde{T}_0 \leq 1\}=\{V(B)_1 \leq 0\}$, where $\widetilde{T}_0:=\inf\{t>0; V(B)_t<0\}$. Formally, it means that we retrieve the information at time $1$ from some prior time, which violates the Markov property.
\subsection{$V(B)$ is semi-martingale -- a conceptual approach}
\label{s42}
When a process is Markov with state space $\mathbb{R}^d$, sufficient and necessary conditions for it to be semi-martingale are given by Cinlar et al. \cite{CJSP}. But $V(B)$ is not Markov. Thus, whether it is semi-martingale cannot be judged by classical Markov-semi-martingale procedures. Here we provide a soft argument to show that $V(B)$ is indeed semi-martingale using Denisov's decomposition of Brownian motion together with Bichteler-Dellacherie's characterization of semi-martingales.

We recall a path decomposition of Brownian motion, which characterizes the Vervaat transform. Let $A$ be the almost sure arcsine split such that $1-A:=\argmin_{t \in [0,1]} B_t$. 
The following theorem is due to Denisov \cite{Denisov}:
\begin{theorem}\cite{Denisov}
\label{Denisov} 
 Given $A$, which is arcsine distributed, i.e. $f_A(a)=\frac{1}{\pi \sqrt{a(1-a)}}$, the path is decomposed into two (conditionally) independent pieces:
\begin{itemize}
\item
$(B_{1-A-t}-B_{1-A}; 0 \leq t \leq 1-A)$  is Brownian meander of length $1-A$;
\item
$(B_t-B_{1-A}; 1-A \leq t \leq 1)$ is Brownian meander of length $A$.
\end{itemize}
\end{theorem}

We turn to some results of the classical semi-martingale theory. 
Given a filtration $(\mathcal{F}_t)_{t \geq 0}$, a process $H$ is said to be simple predictable if $H$ has a representation
$$\forall t \in [0,1],~H_t=H_01_{\{0\}}(t)+\sum_{i=1}^{n-1} H_i 1_{(t_i,t_{i+1}]}(t)$$
where $H_i \in \mathcal{F}_{t_i}$ and $|H_i|< \infty$ a.s. for $0=t_1 \leq t_2 \leq...\leq t_n \leq \infty$. 
Denote $\mathcal{S}$ the collection of simple predictable processes and $\mathcal{B}=\{H \in \mathcal{S}: |H| \leq 1\}$. For a given process, we define a linear mapping $I_X: \mathcal{S} \rightarrow \mathbb{L}^0$ by
$$I_X(H)=H_0X_0+\sum_{i=1}^{n-1}H_i(X_{t_i}-X_{t_{i-1}}),$$
for $H\in \mathcal{S}$. The following theorem, proved independently by Bichteler \cite{Bich} and Dellacherie \cite{Dellacherie}, provides a useful characterization of semi-martingales. We refer readers to Protter \cite[Chapter III]{Protter}, and Dellacherie and Meyer \cite[Section $4$, Chapter VIII]{DM2} for further discussions, and to Beiglb{\"o}ck et al. \cite{BS,BSV} for short proofs.
\begin{theorem} \cite{Bich,Dellacherie}
\label{BD}
An adapted, c\`{a}dl\`{a}g process $X$ is semi-martingale if and only if $I_X(\mathcal{B})$ is bounded in probability, that is 
$$\lim_{\eta \rightarrow \infty} \sup_{H \in \mathcal{B}} \mathbb{P}(|I_X(H)| \geq \eta)=0.$$
\end{theorem}
\begin{proposition}
\label{semi}
$(V(B)_t; 0 \leq t \leq 1)$ is semi-martingale.
\end{proposition}
\begin{proof} Fix $H \in \mathcal{B}$ and $\eta>0$,
\begin{equation}
\label{19}
\mathbb{P}(|I_{V(B)}(H)|>\eta)=\int_0^1 \mathbb{P}(|I_{V(B)|A=a}(H)|>\eta) \frac{1}{\pi \sqrt{a(1-a)}}da.
\end{equation}
Note that $(V(B)|A=1)$ is Brownian meander of length $1$. It has the same distribution as the three dimensional Bessel bridge from $0$ to $\rho$, which is Rayleigh distributed $\mathbb{P}(\rho \in dx)=x \exp(-\frac{x^2}{2})dx$. According to Imhof \cite{Imhof}, and Az\'{e}ma and Yor \cite{AY}, $(V(B)|A=1)$ is semi-martingale and so is $(V(B)|A=0)$. By Theorem \ref{BD}, 
$$\lim_{\eta \rightarrow \infty} \sup_{H \in \mathcal{B}} \mathbb{P}(|I_{V(B)|A=1}(H)| \geq \eta)=0 \quad \mbox{and} \quad \lim_{\eta \rightarrow \infty} \sup_{H \in \mathcal{B}} \mathbb{P}(|I_{V(B)|A=0}(H)| \geq \eta)=0.$$
From \eqref{19}, to prove $\sup_{H \in \mathcal{B}} \mathbb{P}(|I_{V(B})(H)| \geq \eta) \rightarrow 0$, we need some uniform control on $\sup_{H \in \mathcal{B}} \mathbb{P}(|I_{V(B)|A=a}(H)|>\eta)$ for all $a \in [0,1]$.
\begin{lemma}
\label{control}
Let $a \in (0,1]$. Then
$$\sup_{H \in \mathcal{B}} \mathbb{P}(|I_{V(B)|A=a}(H)| > \eta) \leq  \sup_{H \in \mathcal{B}} \mathbb{P}\left(|I_{V(B)|A=1}(H)| > \frac{\eta}{2}\right)+\sup_{H \in \mathcal{B}} \mathbb{P}\left(|I_{V(B)|A=0}(H)| > \frac{\eta}{2}\right).$$
\end{lemma}
\begin{proof} Observe that $I_{V(B)|A=a}(H)=I_1^a+I_2^a$, where
$$I_1^a:=\sum_i H_i(V(B)_{\tau_{i+1} \wedge a}-V(B)_{\tau_{i+1} \wedge a}) \quad \mbox{and} \quad I_2^a:=\sum_i H_i(V(B)_{\tau_{i+1} \vee a}-V(B)_{\tau_{i+1} \vee a}).$$
Then we have,
$$\mathbb{P}(|I_{V(B)|A=a}(H)|>\eta) \leq \mathbb{P}\left(|I_1^a| >\frac{\eta}{2}\right)+\mathbb{P}\left(|I_2^a| >\frac{\eta}{2}\right).$$
Let $\widetilde{I_1^a}:=\frac{I_1^a}{\sqrt{a}}$, and $(B^{me}_t; 0\leq t \leq 1)$ be Brownian meander of length $1$. By Theorem \ref{Denisov}, 
$$
\widetilde{I_1^a}=\sum_i H_i \frac{V(B)_{\tau_{i+1} \wedge a} -V(B)_{\tau_i \wedge a}}{\sqrt{a}}=\sum_i H_i (B^{me}_{\widetilde{\tau}_{i+1}}-B^{me}_{\widetilde{\tau}_i}), 
$$
where $\widetilde{\tau}_i:=\frac{\tau_{i} \wedge a}{a}$, and $H_i$ is $\mathcal{F}^{me}_{\widetilde{\tau}_i}$-adapted.
Consequently,
$$\mathbb{P}\left(|I_1^a| > \frac{\eta}{2}\right) \leq \mathbb{P}\left(|\widetilde{I_1^a}| >\frac{\eta}{2}\right) \leq \sup_{H \in \mathcal{B}} \mathbb{P}\left(|I_{V(B)|A=1}(H)| > \frac{\eta}{2}\right).$$
Similarly, we obtain:
$$\mathbb{P}\left(|I_2^a| > \frac{\eta}{2}\right)  \leq \sup_{H \in \mathcal{B}} \mathbb{P}\left(|I_{V(B)|A=0}(H)| > \frac{\eta}{2}\right).$$ 
\end{proof}
It is straightforward that Proposition \ref{semi} follows from \eqref{19} and Lemma \ref{control}.  
\end{proof}

As a result, Vervaat bridges are also semi-martingales. This provides an alternative proof of the semi-martingale property of Vervaat bridges obtained in Subsection \ref{s33}.
\begin{corollary}
For each fixed $\lambda \in \mathbb{R}$, $(V(B^{\lambda,br})_t; 0 \leq t \leq 1)$ is semi-martingale.
\end{corollary}
\begin{proof} Fix $H \in \mathcal{B}$ and $\eta>0$,
\begin{equation}
\label{20}
\mathbb{P}(I_{V(B)} (H)>\eta)=\int_{\mathbb{R}}\mathbb{P}(I_{V(B^{\lambda,br})}(H) > \eta)\frac{1}{\sqrt{2 \pi}} \exp\left(-\frac{\lambda^2}{2}\right)d\lambda.
\end{equation}
Note that $V(B^{0,br})$ is Brownian excursion, and is semi-martingale. It suffices to prove that $V(B^{\lambda,br})$ where $\lambda \neq 0$ is semi-martingale. If not the case, then there exists $\epsilon>0$ such that for all $K>0$, we can find $\eta>K$ satisfying 
$$\sup_{H \in \mathcal{B}}\mathbb{P}(I_{V(B^{\lambda,br})}(H) > \eta)>\epsilon.$$
It is not hard to see that $(H, \lambda) \rightarrow \mathbb{P}(I_{V(B^{\lambda,br})}(H) > \eta)$ is jointly continuous in $\mathcal{B} \times (\mathbb{R}\setminus \{0\})$, the full detail of which is left to careful readers. Thus, there exists $H_{\lambda,\epsilon} \in \mathcal{B}$ and $\theta \in (0,|\lambda|)$ such that for all $\bar{\lambda} \in (\lambda-\theta,\lambda+\theta)$,
\begin{equation}
\label{21}
\mathbb{P}(I_{V(B^{\bar{\lambda},br})}(H) > \eta)>\frac{\epsilon}{2}.
\end{equation}
Injecting \eqref{21} into \eqref{20}, we obtain:
$$\mathbb{P}(I_{V(B)}(H) > \eta)>\frac{\epsilon}{2}\int_{\lambda-\theta}^{\lambda+\theta}\frac{1}{\sqrt{2 \pi}} \exp\left(-\frac{\lambda^2}{2}\right)d\lambda,$$
which violates that fact that $V(B)$ is a semi-martingale.  
\end{proof}

Note that one can hardly derive the explicit decomposition by Bichteler-Dellacherie's approach. 
A crucial step of the method is to find $\mathbb{Q}$, which is equivalent to $\mathbb{P}$, such that $X$ is $\mathbb{Q}-$quasi-martingale, see e.g. Protter \cite[Chapter III]{Protter}. 
By Rao's theorem, $X$ is $\mathbb{Q}-$semi-martingale, which is also $\mathbb{P}-$semi-martingale. However, Rao's theorem relies on Doob-Meyer's decomposition, which in general does not give explicit expressions for two decomposed terms. 
Nevertheless, we provide the canonical decomposition of $V(B)$ in the next section.
\subsection{Semi-martingale decomposition of $V(B)$}
\label{s43}
In this part, we use extensively the notations defined in Subsection \ref{s33}. Since $V(B)_{1}=B_{1}$ almost surely, there exists $\varepsilon>0$ such that for all $t\in(0,\varepsilon)$, $V(B)_{t}>0$. Let
\begin{displaymath}
\widetilde{T}_{0}:=\inf\lbrace t\in (0,1]; V(B)_{t}=0\rbrace .
\end{displaymath}
Then $\mathbb{P}(\widetilde{T}_{0}\leq 1)=\dfrac{1}{2}$, and
$\lbrace \widetilde{T}_{0}\leq 1\rbrace = \lbrace V(B)_{1}\leq 0\rbrace$. Conditional on $\widetilde{T}_{0}\leq 1$, $\widetilde{T}_{0}$ follows the arcsine distribution $1_{0< t<1}\dfrac{dt}{\pi\sqrt{t(1-t)}}$. Conditionally on $\widetilde{T}_{0}\leq 1$ and the value of $\widetilde{T}_{0}$, $(V(B)_{t};0\leq t\leq \widetilde{T}_{0})$
has the distribution $\mathbb{Q}^{\widetilde{T}_{0}}_{0,0}$, and is independent of $(V(B)_{t};\widetilde{T}_{0}\leq t\leq 1)$.
The joint distribution of  $(V(B)_{1}, \widetilde{T}_{0})$ on the event $\widetilde{T}_{0}\leq 1$ is given by
\begin{displaymath}
\dfrac{1_{\lambda <0}\,d \lambda}{\sqrt{2\pi}}\exp\left(-\dfrac{\lambda^{2}}{2}\right)
\dfrac{\vert\lambda\vert}{\sqrt{2\pi t(1-t)^{3}}}\exp\left(-\dfrac{\lambda^{2}t}{2(1-t)}\right)\,1_{0<t<1}dt.
\end{displaymath}
Thus, the distribution of $V(B)_{1}$ conditionally on $\widetilde{T}_{0}=\tilde{t}_{0}$ is:
\begin{equation}
\label{CondDensEndpoint}
\dfrac{\vert\lambda\vert}{1-\tilde{t}_{0}}
\exp\left(-\dfrac{\lambda^{2}}{2(1-\tilde{t}_{0})}\right)
1_{\lambda <0} d\lambda.
\end{equation}

For the canonical decomposition of $V(B)$, we split the task into two: the decomposition of $(V(B)_{t};0\leq t\leq \widetilde{T}_{0})$ and that of 
$(V(B)_{t};\widetilde{T}_{0}\leq t\leq 1)$. We start with the latter. Let $(\widetilde{M}_{t};t \geq 0)$ be the process defined as
\begin{displaymath}
\widetilde{M}_{t}:=\min_{[0,t]}V(B).
\end{displaymath}

\begin{lemma}
\label{LemEqQueue}
Let $V:=V(B)$. Conditionally on the value of $\widetilde{T}_{0}$,
\begin{displaymath}
\left(V_{t}+\int_{\widetilde{T}_{0}}^{t}\dfrac{V_{s}-\widetilde{M}_{s}}{1-s}\,ds\right)
_{\widetilde{T}_{0}\leq t\leq 1}
\end{displaymath}
is Brownian motion
\end{lemma}
\begin{proof}
Let $\tilde{t}_{0}>0$. Let $(B'_{t})_{t\geq 0}$ be Brownian motion starting from $0$, and
\begin{displaymath}
M'_{t}:=\min_{[0,t]}B'.
\end{displaymath}
For any $t\in [\tilde{t}_{0},1)$, the distribution of $(V(B)_{s};\tilde{t}_{0}\leq s\leq t)$ conditioned on $\widetilde{T}_{0}=\tilde{t}_{0}$, is absolutely continuous with respect to that of $(B'_{s};0\leq s\leq t-\tilde{t}_{0})$. The corresponding density is:
\begin{equation}
\label{DensityReste}
\int_{-\infty}^{M'_{t-\tilde{t}_{0}}}\dfrac{g_{1-t}(B'_{t-\tilde{t}_{0}}-\lambda)}
{g_{1-\tilde{t}_{0}}(|\lambda|)}
\dfrac{\vert\lambda\vert}{1-\tilde{t}_{0}}
\exp\left(-\dfrac{\lambda^{2}}{2(1-\tilde{t}_{0})}\right)d\lambda.
\end{equation}
where $g_t$ is defined as in \eqref{Eq4}. Then the expression \eqref{DensityReste} can be simplified as
$$
\sqrt{\dfrac{1-\tilde{t}_{0}}{(1-t)^{3}}}
\int_{-\infty}^{M'_{t-\tilde{t}_{0}}}(B'_{t-\tilde{t}_{0}}-\lambda)
\exp\left(-\dfrac{(B'_{t-\tilde{t}_{0}}-\lambda)^{2}}
{2(1-t)}\right) d\lambda
=\sqrt{\dfrac{1-\tilde{t}_{0}}{1-t}}
\exp\left(-\dfrac{(B'_{t-\tilde{t}_{0}}-
M'_{t-\tilde{t}_{0}})^{2}}{2(1-t)}\right).
$$
Applying Girsanov's theorem, we get the result of the lemma.  
\end{proof}

Next we deal with the canonical decomposition of $(V(B)_{t\wedge\widetilde{T}_{0}};0\leq t\leq 1)$.
As an auxiliary problem, we study the canonical decomposition of $(\xi_{t};t\geq 0)$ defined as follows. With probability $\frac{1}{2}$, $\xi$ is the three dimensional Bessel process starting from $0$. For $t\in (0,1)$, with infinitesimal probability $\frac{dt}{2\pi\sqrt{t(1-t)}}$, $\xi$ is positive excursion of length $t$, absorbed at $0$ after time $t$. For any $t\in(0,1)$, the distribution of $(V(B)_{s\wedge\widetilde{T}_{0}};0\leq s\leq t)$ is absolutely continuous with respect to that of $(\xi_{s};0\leq s\leq t)$. The following lemma is a variant of Proposition \ref{keylem}.
\begin{proposition}
\label{LemXi}
Let
$
T^{\xi}_{0}:=\inf\lbrace t>0\vert \xi_{t}=0\rbrace
$, and
\begin{displaymath}
J_{t}(y):=\int_{t\leq s\leq 1}\dfrac{ds}{\pi\sqrt{s(1-s)}}\dfrac{\tilde{q}_{s-t}(0,y)}{\tilde{q}_{s}(0,0)},
\end{displaymath}
\begin{displaymath}
\mathring{J}_{t}(y):=\int_{t\leq s\leq 1}\dfrac{ds}{\pi(s-t)\sqrt{s(1-s)}}\dfrac{\tilde{q}_{s-t}(0,y)}
{\tilde{q}_{s}(0,0)}.
\end{displaymath}
The process
\begin{displaymath}
\label{StoppedBM}
(Y_{t})_{t\geq 0}:=\left(\xi_{t}-\int_{0}^{t\wedge T^{\xi}_{0}}\dfrac{ds}{\xi_{s}}+
\int_{0}^{t\wedge T^{\xi}_{0}}\dfrac{\xi_{s}\mathring{J}_{s}(\xi_{s})}{1+J_{s}(\xi_{s})}
\,ds\right)_{t\geq 0}
\end{displaymath}
is Brownian motion with respect to the filtration of $\xi$, stopped at time $T^{\xi}_{0}$.
\end{proposition}
\begin{proof}
Let $\varepsilon\in(0,1)$. Introduce $(B^{\varepsilon}_{t};t\geq 0)$ Brownian motion with the starting point $B^{\varepsilon}_{0}\stackrel{(d)}{=} \xi_{\varepsilon\wedge T^{\xi}_{0}}$. The density of this distribution on $(0,+\infty)$ (total mass $<1$) is:
\begin{displaymath}
\mu_{\varepsilon}(x)=\dfrac{\tilde{q}_{\varepsilon}(0,x)x^{2}}{2}
\left(1+J_{\varepsilon}(x)\right).
\end{displaymath}
Let $T^{\varepsilon}_{0}$ be the first time $B^{\varepsilon}$ hits $0$. For any $\varepsilon \leq t \leq 1$, the distribution of $(\xi_{s};\varepsilon\leq s\leq t)$ is absolutely continuous with respect to that of $(B^{\varepsilon}_{(s-\varepsilon)\wedge T^{\varepsilon}_{0}};\varepsilon\leq s\leq t)$. The density $\mathfrak{D}^{\varepsilon}_{t}$ is:
\begin{align*}
&\quad 1_{B^{\varepsilon}_{0}=0}+
\dfrac{1_{T^{\varepsilon}_{0}\leq t-\varepsilon, B^{\varepsilon}_{0}>0} \cdot  \tilde{q}_{\varepsilon}(0,B^{\varepsilon}_{0})
B^{\varepsilon 2}_{0}\tilde{q}_{T^{\varepsilon}_{0}}(0,B^{\varepsilon}_{0})}
{\tilde{q}_{T^{\varepsilon}_{0}+\varepsilon}(0,0)2\pi\sqrt{T^{\varepsilon}_{0}(1-T^{\varepsilon}_{0})}
\mu_{\varepsilon}(B^{\varepsilon}_{0}) g_{T^{\varepsilon}_{0}}(B^{\varepsilon}_{0})}+\dfrac{1_{T^{\varepsilon}_{0}>t-\varepsilon}}{\mu_{\varepsilon}(B^{\varepsilon}_{0})B^{\varepsilon}_{0}B^{\varepsilon}_{t-\varepsilon}\tilde{q}_{t-\varepsilon}(B^{\varepsilon}_{0},B^{\varepsilon}_{t-\varepsilon})}\\
&\quad \quad \quad \quad \quad \times
\dfrac{\tilde{q}_{\varepsilon}(0,B^{\varepsilon}_{0})B^{\varepsilon 2}_{0}
\tilde{q}_{t-\varepsilon}(B^{\varepsilon}_{0},B^{\varepsilon}_{t-\varepsilon})B^{\varepsilon 2}_{t-\varepsilon}}{2}
\times \left(1+\int_{t}^{1}\dfrac{ds}{\pi\sqrt{s(1-s)}}\dfrac{\tilde{q}_{s-t}
(0,B^{\varepsilon}_{t-\varepsilon})}{\tilde{q}_{s}(0,0)}\right)\\
&=1_{B^{\varepsilon}_{0}=0}+\dfrac{1_{T^{\varepsilon}_{0}\leq t-\varepsilon, B^{\varepsilon}_{0}>0} \cdot \tilde{q}_{\varepsilon}(0,B^{\varepsilon}_{0})
B^{\varepsilon}_{0}}
{\mu_{\varepsilon}(B^{\varepsilon}_{0})\pi\sqrt{T^{\varepsilon}_{0}(1-T^{\varepsilon}_{0})}
\tilde{q}_{T^{\varepsilon}_{0}+\varepsilon}(0,0)}+1_{T^{\varepsilon}_{0}>t-\varepsilon}\dfrac{\tilde{q}_{\varepsilon}(0,B^{\varepsilon}_{0})B^{\varepsilon}_{0}
B^{\varepsilon}_{t-\varepsilon}}{2\mu_{\varepsilon}(B^{\varepsilon}_{0})}\times
\left(1+J_{t}(B^{\varepsilon}_{t-\varepsilon})\right).
\end{align*}
Note that $(\mathfrak{D}^{\varepsilon}_{t};t\geq 0)$ is continuous, and there is no discontinuity at $T^{\varepsilon}_{0}+\varepsilon$. This follows from the fact that as $y \rightarrow 0$, the convolution kernel
$
\dfrac{y}{2}\tilde{q}_{u}(0,y)\,1_{u>0}\,du
$
is an approximation to the delta function.
Moreover,
$
\dfrac{\partial J_{t}(y)}{\partial y}=-y\mathring{J}_{t}(y)
$
for $t\in(0,1)$. 
We get the quadratic covariation
\begin{displaymath}
d\left[\log(\mathfrak{D}^{\varepsilon}_{t\wedge 
T_0^{\varepsilon}}),
B_{(t-\varepsilon)\wedge T_0^{\varepsilon}}
^{\varepsilon}\right]=1_{T_0^{\varepsilon} > t-\varepsilon}
\dfrac{\tilde{q}_{\varepsilon}(0,B^{\varepsilon}_{0})B^{\varepsilon}_{0}
(1+J_t(B^{\varepsilon}_{t-\varepsilon})-
(B^{\varepsilon}_{t-\varepsilon})^{2} \mathring{J}_t(B^{\varepsilon}_{t-\varepsilon}))}
{2\mu_{\varepsilon}(B^{\varepsilon}_{0})}.
\end{displaymath}
Applying Girsanov's theorem, we obtain that $(Y_{t};t\geq \varepsilon)$ is a continuous martingale relative to the filtration of $(\xi_{t};t\geq \varepsilon)$ with quadratic variation $(t-\varepsilon)\wedge
(T^{\xi}_{0}-\varepsilon)^{+}$. Since this holds for all $\varepsilon$ sufficiently small, this implies the result.  
\end{proof}

We introduce the functionals $\Phi(t,\gamma)$ and $\mathring{\Phi}(t,\gamma)$ where $t$ is a time and $\gamma$ a continuous path:
\begin{displaymath}
\Phi(t,\gamma):=\dfrac{2}{\sqrt{2\pi}}\int_{0}^{+\infty}
\Phi^{\lambda}(t,\gamma(t),\sup\lbrace s\in [0,t]\vert \gamma(s)\leq \lambda\rbrace)
\exp\left(-\dfrac{\lambda^{2}}{2}\right)\,d\lambda,
\end{displaymath}
\begin{displaymath}
\mathring{\Phi}(t,\gamma):=\dfrac{2}{\sqrt{2\pi}}\int_{0}^{+\infty}
\partial_{2}\Phi^{\lambda}(t,\gamma(t),\sup\lbrace s\in [0,t]\vert \gamma(s)\leq \lambda\rbrace)
\exp\left(-\dfrac{\lambda^{2}}{2}\right)\,d\lambda,
\end{displaymath}
where $\Phi^{\lambda}$ is defined as \eqref{phi}. For any $t\in(0,1)$, the distribution of $(V(B)_{s\wedge\widetilde{T}_{0}};0\leq s\leq t)$ is absolutely continuous with respect to that of $(\xi_{s};0\leq s\leq t)$ with density
\begin{equation}
\label{FinalDensity}
D_{t}=1_{T^{\xi}_{0}\leq t}+\dfrac{\Phi(t,\xi)+J_{t}(\xi_{t})}{1+J_{t}(\xi_{t})}\,1_{T^{\xi}_{0}> t}.
\end{equation}

The following lemma provides estimates on $\Phi^{\lambda}$ and its derivatives $\partial_1\Phi^{\lambda},\partial_2\Phi^{\lambda}$.
\begin{lemma}
\label{LemBound}
There are positive functions $c_1(t)$ and $c_2(t)$ bounded on intervals of the form $[0,1-\varepsilon]$, such that for all 
$\lambda>0$, $y>0$, $\theta\leq t\in [0,1)$:
\begin{displaymath}
\Phi^{\lambda}(t,y,\theta)\leq c_1(t)\exp\left(\frac{\lambda^{2}}{2}\right)
\exp\left(-\dfrac{(y-\lambda)^{2}}{2(1-t)}\right),
\end{displaymath} 
and
\begin{displaymath}
|\partial_2 \Phi^{\lambda}(t,y,\theta)|\leq c_2(t)(1+\lambda^2)\left(y+\frac{1}{y}\right)\exp\left(\frac{\lambda^{2}}{2}\right)
\exp\left(-\dfrac{(y-\lambda)^{2}}{2(1-t)}\right).
\end{displaymath} 
\end{lemma}
\begin{proof} Note that $\Phi^{\lambda}(t,y,\theta)\leq \Phi^{1,\lambda}(t,y)+1_{y>\lambda}\Phi^{2,\lambda}(t,y)$. For $y>\lambda$,
\begin{displaymath}
\Phi^{2,\lambda}(t,y)\leq \dfrac{1}{\sqrt{(1-t)^{3}}}\exp\left(\frac{\lambda^{2}}{2}\right)
\exp\left(-\dfrac{(y-\lambda)^{2}}{2(1-t)}\right)
\end{displaymath}
and
\begin{equation*}
\begin{split}
\Phi^{1,\lambda}(t,y)\leq&\dfrac{1}{2\sqrt{2}y}\exp\left(\frac{\lambda^{2}}{2}\right)
\exp\left(-\dfrac{(y-\lambda)^{2}}{2(1-t)}\right)
\int_{\frac{(y-\lambda)^{2}}{2(1-t)}}^{\frac{(y+\lambda)^{2}}{2(1-t)}}\dfrac{du}{\sqrt{u}}\\
=&\dfrac{1}{\sqrt{1-t}}\dfrac{\min(y,\lambda)}{y}
\exp\left(\frac{\lambda^{2}}{2}\right)
\exp\left(-\dfrac{(y-\lambda)^{2}}{2(1-t)}\right). 
\end{split}
\end{equation*}
In addition, for $y>\lambda$ we obtain,
\begin{align}
|\partial_2 \Phi^{2,\lambda}(t,y)| &=\frac{1}{\sqrt{(1-t)^3}} \left| \frac{\lambda}{y^2}-\frac{(y-\lambda)^2}{(1-t)y} \right| \exp\left(\frac{\lambda^2}{2}\right) \exp\left(-\frac{(y-\lambda)^2}{2(1-t)}\right)  \notag\\
                                                         & \leq \frac{1}{\sqrt{(1-t)^3}} \left(\frac{1}{y}+\frac{y^2+\lambda^2}{(1-t)y} \right)\exp\left(\frac{\lambda^2}{2}\right) \exp\left(-\frac{(y-\lambda)^2}{2(1-t)}\right)  \notag
\end{align}
and
\begin{multline*}
|\partial_2 \Phi^{1,\lambda}(t,y)| \leq \frac{1}{2 \sqrt{2}} \exp\left(\frac{\lambda^2}{2}\right) \Bigg( \frac{1}{y^2}  \int_{\frac{(y-\lambda)^{2}}{2(1-t)}}^{\frac{(y+\lambda)^{2}}{2(1-t)}}\dfrac{e^{-u}}{\sqrt{u}}du 
\\+ \frac{\sqrt{2}}{y \sqrt{1-t}} \left| \exp\left(-\frac{(y-\lambda)^2}{2(1-t)} \right)- \exp\left(-\frac{(y-\lambda)^2}{2(1-t)} \right) \right| \Bigg) 
\end{multline*}
$$\leq \frac{1+2 \sqrt{2}}{2 \sqrt{2} y \sqrt{1-t}} \exp\left(\frac{\lambda^2}{2}\right) \exp\left(-\frac{(y-\lambda)^2}{2(1-t)}\right), \quad\quad\quad\quad\quad\quad$$
which permits to have the desired estimation.  
\end{proof}

\begin{lemma}
\label{LemQuadVar}
For $t\in[0,1)$,
\begin{displaymath}
[\Phi(\cdot ,\xi),\xi]_{t}=\int_{0}^{t\wedge T^{\xi}_{0}}\mathring{\Phi}(s,\xi)\,ds.
\end{displaymath}
\end{lemma}
\begin{proof} It is clear that the quadratic variation $[\Phi(\cdot ,\xi),\xi]_{t}$ does not increase for $t\geq T^{\xi}_{0}$. We only need to show that for a three dimensional Bessel process $(R_{t};t\geq 0)$,
\begin{equation}
\label{EqQuadVar}
[\Phi(\cdot ,R),R]_{t}=\int_{0}^{t}\mathring{\Phi}(s,R)\,ds.
\end{equation}
Indeed, given any $t<T<1$, the distribution of $(\xi_{s};0\leq s\leq t)$ on the event $T^{\xi}_{0}>T$ is absolutely continuous with respect to that of $(R_{s};0\leq s\leq t)$.
For any $\lambda>0$ $(\Phi^{\lambda}(t,R_{t},\theta^{\lambda}_{t});0\leq t<1)$ is a positive martingale with mean $1$. Applying Fubini's theorem, we obtain that $(\Phi(t,R);0\leq t<1)$ is a positive martingale with mean $1$. Let 
$(W_{t};t\geq 0)$ be the Brownian motion martingale part of $(R_{t};t\geq 0)$. To prove \eqref{EqQuadVar} we need to show that the process
\begin{equation}
\label{ShowTrueMart}
\left(\Phi(t,R)W_{t}-\int_{0}^{t}\mathring{\Phi}(s,R)\,ds\right)_{0\leq t<1}
\end{equation}
is a martingale. Lemma \ref{LemMartDec} ensures that for any $\lambda>0$ the process
\begin{equation}
\label{EqMartLoc}
\left(\Phi^{\lambda}(t,R_{t},\theta^{\lambda}_{t})W_{t}-\int_{0}^{t}
\partial_{2}\Phi^{\lambda}(s,R_{s},\theta^{\lambda}_{s})\,ds\right)_{0\leq t<1}
\end{equation}
is a local martingale. Next we show that \eqref{EqMartLoc} is a true martingale and not just a local one. It suffices to bound the expectation of its supremum and dominated convergence theorem permits to conclude. According to Burkholder-Davis-Gundy's inequality, $\exists C>0$ such that 
\begin{align*}
&~\quad \mathbb{E}\Bigg[ \Bigg(\sup_{0 \leq s \leq t} \Bigg|\Phi^{\lambda}(t,R_{t},\theta^{\lambda}_{t})W_{t}-\int_{0}^{t} \partial_{2}\Phi^{\lambda}(s,R_{s},\theta^{\lambda}_{s})\,ds \Bigg| \Bigg)^2\Bigg] \\
&\leq C \mathbb{E} \left[\int_0^t \left( \partial_2 \Phi^{\lambda}(s,R_s,\theta_s^{\lambda})^{2}W_s^{2} +
\Phi^{\lambda}(s,R_s,\theta_s^{\lambda})^{2}\right)ds\right]\\
&=C \int_0^t\left(\mathbb{E}\left[\partial_2 \Phi^{\lambda}(s,R_s,\theta_s^{\lambda})^{2}W_s^{2}\right]+
\mathbb{E}\left[\Phi^{\lambda}(s,R_s,\theta_s^{\lambda})^{2}\right]\right)ds.
\end{align*}
From Lemma \ref{LemMartDec} and the bound in Lemma \ref{LemBound} follows that $(\Phi^{\lambda}(t,R_{t},\theta^{\lambda}_{t});0\leq t<1)$ is a square integrable martingale, and
\begin{equation*}
1+\mathbb{E}\left[\int_{0}^{s}
\partial_{2}\Phi^{\lambda}(u,R_{u},\theta^{\lambda}_{u})^{2}\,du\right]=
\mathbb{E}\left[\Phi^{\lambda}(s,R_{s},\theta^{\lambda}_{s})^{2}\right]\leq
c_1(s)^{2}\exp\left(\lambda^{2}\right)
\mathbb{E}\left[\exp\left(-\dfrac{(R_{s}-\lambda)^{2}}{(1-s)}\right)\right],
\end{equation*}
which is integrable on $(0,t)$ for any $0 \leq t <1$. Moreover, by Cauchy-Schwarz's inequality, 
$$
\mathbb{E}[\partial_2 \Phi^{\lambda}(s,R_s,\theta_s^{\lambda})^2W_s^2] \leq \mathbb{E}[\partial_2 \Phi^{\lambda}(s,R_s,\theta_s^{\lambda})^4]^{\frac{1}{2}} \mathbb{E}[W_s^4]^{\frac{1}{2}}.
$$
The problem of integrability may only occur at $0$. However, by the bound of $\partial_2 \Phi^{\lambda}$ in Lemma \ref{LemBound}, we know that $\mathbb{E}[\partial_2 \Phi^{\lambda}(s,R_s,\theta_s^{\lambda})^4]=\mathcal{O}(\frac{1}{s^2})$ as $s \rightarrow 0$. Thus the above term is also integrable on $(0,t)$ for $0 \leq t \leq 1$. We have proved that \eqref{EqMartLoc} is a martingale for any $\lambda >0$. Again by Cauchy-Schwarz's inequality,
\begin{multline*}
\mathbb{E}\left[\left\vert
\Phi^{\lambda}(t,R_{t},\theta^{\lambda}_{t})W_{t}-\int_{0}^{t}
\partial_{2}\Phi^{\lambda}(s,R_{s},\theta^{\lambda}_{s})\,ds
\right\vert\right]\leq
2\sqrt{t}\mathbb{E}\left[
\int_{0}^{t}\partial_{2}\Phi^{\lambda}(s,R_{s},\theta^{\lambda}_{s})^{2}\,ds
\right]^{\frac{1}{2}}\\
\leq 2\sqrt{t} c_{1}(t)\exp\left(\dfrac{\lambda^{2}}{2}\right)
\mathbb{E}\left[\exp\left(-\dfrac{(R_{t}-\lambda)^{2}}{(1-t)}\right)\right]^{\frac{1}{2}}< \infty.
\end{multline*}
Thus, the expectation of the absolute value of the martingale \eqref{EqMartLoc} is integrable with respect to $\dfrac{2}{\sqrt{2\pi}}\exp\left(-\dfrac{\lambda^{2}}{2}\right)\,1_{\lambda >0}\,d\lambda$. By Fubini's theorem, it follows that \eqref{ShowTrueMart} is a martingale.  
\end{proof}

\begin{theorem}
Let $\tilde{M}_t:=\min_{[0,t]}V(B)$, $\tilde{T}_0:=\inf\{t>0; V(B)_t=0\}$ and $V:=V(B)$. Then 
$$
\Bigg(V_{t}-\int_{0}^{t\wedge\widetilde{T}_{0}}\dfrac{ds}{V_{s}}+\int_{0}^{t\wedge\widetilde{T}_{0}}\dfrac{\mathring{\Phi}(s,V)+
V_{s}\mathring{J}_{s}(V_{s})}{\Phi(s,V)+J_{s}(V_{s})}\,ds
+\int_{\widetilde{T}_{0}}^t \dfrac{V_{s}-\widetilde{M}_{s}}{1-s}\,ds\Bigg)_{0\leq t\leq 1}
$$
is Brownian motion.
\end{theorem}
\begin{proof} Note that $(D_{t})_{0\leq t\leq 1}$ given by \eqref{FinalDensity} is time-continuous. In particular it follows from Lemma \ref{LemBound} that on the event $T^{\xi}_{0}<1$, as $t$ converges to $T^{\xi}_{0}$ from below and $\xi_{t}$ converges to $0$, $\Phi(t,\xi)$ remains bounded. In addition, $J_{t}(\xi_{t})$ tends to $+\infty$ at $T^{\xi}_{0}$. Hence
\begin{displaymath}
\lim_{t\rightarrow T^{\xi}_{0}}\dfrac{\Phi(t,\xi)+J_{t}(\xi_{t})}{1+J_{t}(\xi_{t})}=1,
\end{displaymath}
and $D_{t}$ is continuous at $t=T^{\xi}_{0}$. Using the canonical decomposition of $(\xi_{t};t\geq 0)$ given by Lemma \ref{LemXi}, and applying Girsanov's theorem combined with Lemma \ref{LemQuadVar}, we obtain that
\begin{displaymath}
\left(V(B)_{t\wedge\widetilde{T}_{0}}-\int_{0}^{t\wedge\widetilde{T}_{0}}\dfrac{ds}{V(B)_{s}}+
\int_{0}^{t\wedge\widetilde{T}_{0}}\dfrac{\mathring{\Phi}(s,V(B))+
V(B)_{s}\mathring{J}_{s}(V(B)_{s})}{\Phi(s,V(B))+J_{s}(V(B)_{s})}\,ds\right)_{t\geq 0}
\end{displaymath}
is a martingale with quadratic variation $t\wedge\widetilde{T}_{0}$. Finally, Lemma \ref{LemEqQueue} describes the canonical decomposition of $V(B)$ after the stopping time $\widetilde{T}_{0}$.  
\end{proof}
\subsection{Expectation and variance of $V(B)$}
\label{s44}
In this subsection, we calculate the first two moments of the Vervaat transform of Brownian motion.
\begin{proposition}
$\forall t \in [0,1]$, we have:
\begin{equation}
\label{23}
\mathbb{E}V(B)_t=\sqrt{\frac{8}{\pi}} (\sqrt{t}+\sqrt{1-t}-1);
\end{equation}
\begin{equation}
\label{24}
\mathbb{E}(V(B)_t^2)=3t+\frac{4-8t}{\pi} \arcsin\sqrt{t}-\frac{4}{\pi}\sqrt{t(1-t)}.
\end{equation}
\end{proposition}

The computation is based on Denisov's decomposition of Brownian motion, Theorem \ref{Denisov}, together with the following identities of Brownian meander, the proof of which is reported to the Appendix:
\begin{proposition}
\label{compmean}
Let $(B^{me}_t, t \in [0,1])$ be Brownian meander of length $1$. We have:
\begin{equation}
\label{25}
\mathbb{E}B^{me}_t=\sqrt{\frac{2}{\pi}} (\sqrt{t(1-t)}+ \arcsin \sqrt{t}).
\end{equation}
\begin{equation}
\label{26}
\mathbb{E}(B^{me}_t)^2=3t-t^2.
\end{equation}
\begin{equation}
\label{27}
\mathbb{E}B^{me}_tB^{me}=2 \sqrt{t}.
\end{equation}
\end{proposition}
\subsubsection{Expectation of $V(B)$}
\quad ~~Recall that $A$ is the almost sure arcsine split such that $1-A:=\argmin_{t \in [0,1]} B_t$. We have:
$$\mathbb{E}V(B)_t=\mathbb{E}(V(B)_t 1_{A>t})+ \mathbb{E}(V(B)_t 1_{A \leq t}).$$
\begin{lemma}
\label{pre1}
\begin{equation}
\label{28}
\forall t \in [0,1],~\mathbb{E}(V(B)_t 1_{A>t})=\sqrt{\frac{2}{\pi}}(\sqrt{1-t}+2 \sqrt{t}-t-1).
\end{equation}
\end{lemma}
\begin{proof} Note that
$$
\mathbb{E}(V(B)_t 1_{A>t})=\int_t^1 \frac{\sqrt{a} \mathbb{E}B^{me}(\frac{t}{a})}{\pi \sqrt{a(1-a)}} da \stackrel{\eqref{25}}{=}\frac{\sqrt{2}}{\pi^{\frac{3}{2}}} (\alpha_1+\alpha_2),                               
$$ 
where $\alpha_1:=\sqrt{t} \int_t^1 \frac{\sqrt{a-t}}{a \sqrt{1-a}}da$ and $\alpha_2:=\int_t^1 \frac{\arcsin\sqrt{\frac{t}{a}}}{\sqrt{1-a}}da.$
Using integration by parts, we get:
$$
\alpha_2=\pi \sqrt{1-t} -\sqrt{t} \int_t^1 \frac{\sqrt{1-a}}{a \sqrt{a-t}}da.
$$ 
Therefore,
\begin{align}
\label{29}
\alpha_1+\alpha_2 &=\pi \sqrt{1-t} + 2 \sqrt{t} \int_t^1 \frac{da}{\sqrt{1-a}\sqrt{a-t}} - \sqrt{t}(t+1) \int_t^1 \frac{da}{a \sqrt{1-a}\sqrt{a-t}} \notag\\
                               &=\pi \sqrt{1-t}+2 \pi \sqrt{t}-\pi(t+1).~~\square \notag
\end{align} 
\end{proof}

Observing the duality 
$(1-A,(V(B)_{1-t}-V(B)_1; 0 \leq t \leq 1)) \stackrel{(d)}{=}(A,(V(B)_t; 0 \leq t \leq 1))$,
we obtain the following result as a corollary:
\begin{corollary}
\label{pre2}
\begin{equation}
\label{33}
\forall t \in [0,1],~\mathbb{E}(V(B)_t 1_{A \leq t})=\sqrt{\frac{2}{\pi}} (\sqrt{1-t}+t-1).
\end{equation}
\end{corollary}

It is not hard to see that \eqref{23} follows from  Lemma \ref{pre1} and Corollary \ref{pre2}.
\subsubsection{Variance of $V(B)$} 
\quad~~Similarly, we split $\mathbb{E}V(B)_t^2$ into $\mathbb{E}(V(B)_t^2 1_{A>t})$ and $\mathbb{E}(V(B)_t^2 1_{A \leq t})$. The formula \eqref{24} follows readily from the next lemmas.
\begin{lemma}
\label{pre}
\begin{equation}
\label{36}
\forall t \in [0,1],~\mathbb{E}(V(B)_t^2 1_{A>t})=3t-\frac{6t}{\pi} \arcsin\sqrt{t}-\frac{2}{\pi} \sqrt{t^3(1-t)}.
\end{equation}
\end{lemma}
\begin{proof} Note that
\begin{align*}
\mathbb{E}(V(B)_t^2 1_{A>t})&=\int_t^1 \frac{a \mathbb{E}(B^{me})^2(\frac{t}{a})}{\pi \sqrt{a(1-a)}} da \\
                                              &\stackrel{(*)}{=}\frac{3t}{\pi} \int_t^1 \frac{1}{\sqrt{a(1-a)}}da -\frac{t^2}{\pi} \int_t^1 \frac{1}{\sqrt{a^3(1-a)}}da \\
                                              &=3t-\frac{6t}{\pi} \arcsin\sqrt{t}-\frac{2}{\pi} \sqrt{t^3(1-t)},
\end{align*}
where $(*)$ follows from \eqref{26}.  
\end{proof}

\begin{lemma}
\begin{equation}
\label{37}
\forall t \in [0,1],~\mathbb{E}(V(B)_t^2 1_{A \leq t})=\frac{4-2t}{\pi} (\arcsin \sqrt{t} -\sqrt{t(1-t)}).
\end{equation}
\end{lemma}
\begin{proof} Note that  
\begin{align}
\mathbb{E}(V(B)_t^2 1_{A \leq t}) &=\mathbb{E}((V(B)_t-V(B)_1)^2 1_{A \leq t})+\mathbb{E}(V(B)_1^2 1_{A \leq t}) \notag\\
                                                         &\quad +2\mathbb{E}((V(B)_t-V(B)_1)V(B)_1 1_{A \leq t}). \notag
\end{align}
Denote $\beta_1, \beta_2$ and $\beta_3$ the three terms on the RHS of the above equation. $\beta_1$ can be easily derived from Lemma \ref{pre}  by change of variables:
\begin{equation}
\label{38}
\beta_1 = \frac{2(1-t)}{\pi}\left(3 \arcsin\sqrt{t}-\sqrt{t(1-t)}\right).  
\end{equation}
By Theorem \ref{Denisov}, for $0 \leq a \leq t$, 
\begin{equation}
\label{39}
\beta_2=\int_0^t \mathbb{E} (V(B)_1^2|A=a) \frac{da}{\pi \sqrt{a(1-a)}}=\int_0^t (2-\pi \sqrt{a(1-a)}) \frac{da}{\pi \sqrt{a(1-a)}}=\frac{4 \arcsin \sqrt{t}}{\pi}-t,
\end{equation}      
and 
$$ \beta_3 =\int_0^t (\gamma_1+\gamma_2)\frac{da}{\pi \sqrt{a(1-a)}},$$
where
\begin{align}
\gamma_1&:=\mathbb{E}\left((V(B)_t-V(B)_1)V(B)_a|A=a\right) \notag\\
       &=\sqrt{a(1-a)} \left(\sqrt{\frac{1-t}{1-a}(1-\frac{1-t}{1-a})} +\arcsin \sqrt{\frac{1-t}{1-a}}\right), \notag
\end{align}
and
\begin{align}
 \gamma_2&:=\mathbb{E}\left((V(B)_t-V(B)_1)(V(B)_1-V(B)_a)|A=a\right) \notag\\
       &=-2 \sqrt{(1-t)(1-a)}, \notag
\end{align}
which are derived from Proposition \ref{compmean}. 
Thus,
\begin{align}
\label{42}
\int_0^t \gamma_1 \frac{da}{\pi \sqrt{a(1-a)}} &= \frac{1}{\pi} \int_{1-t}^1  \left(\sqrt{\frac{1-t}{a}(1-\frac{1-t}{a})} +\arcsin \sqrt{\frac{1-t}{a}}\right) da \notag\\
                                                                                 &=\frac{t}{2}+\frac{2t-3}{\pi} \arcsin \sqrt{t} +\frac{3}{\pi} \sqrt{t(1-t)},
\end{align}
and
\begin{equation}
\label{43}
\int_0^t \gamma_2 \frac{da}{\pi \sqrt{a(1-a)}}= -\frac{4}{\pi}\sqrt{t(1-t)}.
\end{equation}
Combining \eqref{38}, \eqref{39}, \eqref{42} and \eqref{43}, we obtain \eqref{37}.  
\end{proof}
\section{Appendix: Computations of Brownian meander}
Let $(B^{me}_t, t \in [0,1])$ be Brownian meander of length $1$. From Chung \cite{Chung}, we derive the marginal distribution of meander:
\begin{equation}
\label{44}
\mathbb{P}(B^{me}_t \in dx)=t^{-\frac{3}{2}}x \exp\left(-\frac{x^2}{2t}\right) \erf\left(\frac{x}{\sqrt{2(1-t)}}\right)dx.
\end{equation}
where $\erf$ is the error function defined as $\erf(x):=\frac{2}{\sqrt{\pi}} \int_{- \infty}^x \exp(-t^2)dt$.
\begin{proof}[Proof of Proposition \ref{compmean}]
$(1)$. We compute $\mathbb{E}B^{me}_t$, which relies on the following well-known identity found in Gradshteyn and Ryzhik \cite{GR}:
\begin{equation}
\label{46}
\forall a>0,~\int_0^{\infty} x^2 \exp(-ax^2) \erf(x)dx =\frac{\sqrt{a}+(a+1)\arcsin\sqrt{\frac{1}{a+1}}}{2 \sqrt{\pi} a^{\frac{3}{2}}(a+1)} .
\end{equation}
By change of variables, we obtain:
\begin{align}
\mathbb{E}B^{me}_t&~~~~\stackrel{(*)}{=}\int_0^{\infty}   t^{-\frac{3}{2}}y^2 \exp\left(-\frac{y^2}{2t}\right) \erf\left(\frac{y}{\sqrt{2(1-t)}}\right) dy \notag\\
                                   &\stackrel{(\#)}{=}\sqrt{8} \left(\frac{1-t}{t}\right)^{\frac{3}{2}} \int_0^{\infty} x^2 \exp\left(-\frac{1-t}{t}x^2\right) \erf(x) dx \notag\\
                                   &\stackrel{(**)}{=}\sqrt{\frac{2}{\pi}} (\sqrt{t(1-t)}+ \arcsin \sqrt{t}),
\end{align}
where $(*)$ (resp. $(**)$) follows from \eqref{44} (resp. \eqref{46}), and $(\#)$ is obtained by taking $x=\frac{y}{\sqrt{2(1-t)}}$.

$(2)$. To calculate $\mathbb{E}(B^{me}_t)^2$, we make use of the following identity: 
\begin{equation}
\label{47}
\forall a>0,~\int_0^{\infty} x^3 \exp(-ax^2) \erf(x)dx =\frac{2+3a}{4a^2(a+1)^{\frac{3}{2}}},
\end{equation}
which is also found in Gradshteyn and Ryzhik \cite{GR}. The rest is the same as in $(1)$.

$(3)$. Finally, the value of $\mathbb{E}B^{me}_tB^{me}_1$, i.e. the formula \eqref{27}, can be easily derived from Imhof's relation \cite{Imhof} between Brownian meander and the three dimensional Bessel process.  
\end{proof}



\begin{thebibliography}{99}

\bibitem{Abr} Josh Abramson, Jim Pitman, Nathan Ross, and Ger{\'o}nimo Uribe~Bravo. \newblock Convex minorants of random walks and {L}\'evy processes. \newblock {\em Electron. Commun. Probab.}, 16:423--434, 2011. \MR{2831081}

\bibitem{Aldous}
David Aldous.
\newblock Brownian excursion conditioned on its local time.
\newblock {\em Electron. Comm. Probab.}, 3:79--90 (electronic), 1998. \MR{1650567}

\bibitem{AP} David Aldous and Jim Pitman. \newblock The standard additive coalescent. \newblock {\em Ann. Probab.}, 26(4):1703--1726, 1998. \MR{1675063}

\bibitem{AFP}
Sami Assaf, Noah Forman, and Jim Pitman.
\newblock The quantile transform of a simple walk.
\newblock 2013.
\newblock \ARXIV{1307.4967}.

\bibitem{AY}
Jacques. Az{\'e}ma and Marc. Yor.
\newblock \'{E}tude d'une martingale remarquable.
\newblock In {\em S\'eminaire de {P}robabilit\'es, {XXIII}}, volume 1372 of
  {\em Lecture Notes in Math.}, pages 88--130. Springer, Berlin, 1989.

\bibitem{BSV} Mathias Beiglb{\"o}ck, Walter Schachermayer, and Bezirgen Veliyev. \newblock A direct proof of the {B}ichteler-{D}ellacherie theorem and  connections to arbitrage. \newblock {\em Ann. Probab.}, 39(6):2424--2440, 2011. \MR{2932672}

\bibitem{BS} Mathias Beiglb{\"o}ck and Pietro Siorpaes. \newblock Riemann-integration and a new proof of the {B}ichteler-{D}ellacherie  theorem. \newblock {\em Stochastic Process. Appl.}, 124(3):1226--1235, 2014. \MR{3148011}

\bibitem{Bertoin} Jean Bertoin. \newblock A fragmentation process connected to {B}rownian motion. \newblock {\em Probab. Theory Related Fields}, 117(2):289--301, 2000. \MR{1771665}

\bibitem{BCP} Jean Bertoin, Lo{\"{\i}}c Chaumont, and Jim Pitman. \newblock Path transformations of first passage bridges. \newblock {\em Electron. Comm. Probab.}, 8:155--166 (electronic), 2003. \MR{2042754}

\bibitem{BCY} Jean Bertoin, Lo\"{i}c Chaumont, and Marc Yor. \newblock Two chain-transformations and their applications to quantiles. \newblock {\em J. Appl. Probab.}, 34(4):882--897, 1997. \MR{1484022}

\bibitem{Biane} Philippe Biane. \newblock Relations entre pont et excursion du mouvement brownien r\'eel. \newblock {\em Ann. Inst. H. Poincar\'e Probab. Statist.}, 22(1):1--7, 1986. \MR{0838369}

\bibitem{BY} Philippe Biane and Marc Yor. \newblock Quelques pr\'ecisions sur le m\'eandre brownien. \newblock {\em Bull. Sci. Math. (2)}, 112(1):101--109, 1988. \MR{0942801}

\bibitem{Bich} Klaus Bichteler. \newblock Stochastic integration and {$L^{p}$}-theory of semimartingales. \newblock {\em Ann. Probab.}, 9(1):49--89, 1981. \MR{0606798}

\bibitem{Bill} Patrick Billingsley. \newblock {\em Probability and measure}. \newblock Wiley Series in Probability and Mathematical Statistics. John Wiley  \& Sons, Inc., New York, third edition, 1995. \newblock A Wiley-Interscience Publication. \MR{1324786}

\bibitem{Bill2} Patrick Billingsley. \newblock {\em Convergence of probability measures}. \newblock Wiley Series in Probability and Statistics: Probability and  Statistics. John Wiley \& Sons Inc., New York, second edition, 1999. \newblock A Wiley-Interscience Publication. \MR{1700749}

\bibitem{CJ} Philippe Chassaing and Svante Janson. \newblock A {V}ervaat-like path transformation for the reflected {B}rownian  bridge conditioned on its local time at 0. \newblock {\em Ann. Probab.}, 29(4):1755--1779, 2001. \MR{1880241}

\bibitem{Chaumont1999} Lo\"{i}c Chaumont. \newblock A path transformation and its applications to fluctuation theory. \newblock {\em J. London Math. Soc. (2)}, 59(2):729--741, 1999. \MR{1709677}

\bibitem{Chaumont} Lo\"{i}c Chaumont. \newblock An extension of {V}ervaat's transformation and its consequences. \newblock {\em J. Theoret. Probab.}, 13(1):259--277, 2000. \MR{1744984}

\bibitem{CU}
Lo\"{i}c Chaumont and Ger{\'o}nimo Uribe~Bravo.
\newblock Shifting processes with cyclically exchangeable increments at random.
\newblock 2014.
\newblock \ARXIV{1405.1335}.

\bibitem{Chung} Kai~Lai Chung. \newblock Excursions in {B}rownian motion. \newblock {\em Ark. Mat.}, 14(2):155--177, 1976. \MR{0467948}

\bibitem{CJSP}
Erhan {\c{C}}inlar, Jean Jacod, Philip Protter, and Michael Sharpe.
\newblock Semimartingales and {M}arkov processes.
\newblock {\em Z. Wahrsch. Verw. Gebiete}, 54(2):161--219, 1980.

\bibitem{Dassios} Angelos Dassios. \newblock The distribution of the quantile of a {B}rownian motion with drift  and the pricing of related path-dependent options. \newblock {\em Ann. Appl. Probab.}, 5(2):389--398, 1995. \MR{1336875}

\bibitem{Dassiosbis} Angelos Dassios. \newblock Sample quantiles of stochastic processes with stationary and  independent increments. \newblock {\em Ann. Appl. Probab.}, 6(3):1041--1043, 1996. \MR{1410129}

\bibitem{Dassiostri} Angelos Dassios. \newblock On the quantiles of {B}rownian motion and their hitting times. \newblock {\em Bernoulli}, 11(1):29--36, 2005. \MR{2121453}

\bibitem{Dellacherie} Claude Dellacherie. \newblock Un survol de la th\'eorie de l'int\'egrale stochastique. \newblock {\em Stochastic Process. Appl.}, 10(2):115--144, 1980. \MR{0587420}

\bibitem{DM2} Claude Dellacherie and Paul-Andr{\'e} Meyer. \newblock {\em Probabilities and potential. {B}}, volume~72 of {\em  North-Holland Mathematics Studies}. \newblock North-Holland Publishing Co., Amsterdam, 1982. \newblock Theory of martingales, Translated from the French by J. P. Wilson. \MR{0745449}

\bibitem{Denisov} I.~V. Denisov. \newblock A random walk and a {W}iener process near a maximum. \newblock {\em Theory of Probability \& Its Applications}, 28(4):821--824,  1984. \MR{0726906}

\bibitem{ERY} Paul. Embrechts, L.~C.~G. Rogers, and Marc Yor. \newblock A proof of {D}assios' representation of the {$\alpha$}-quantile of  {B}rownian motion with drift. \newblock {\em Ann. Appl. Probab.}, 5(3):757--767, 1995. \MR{1359828}

\bibitem{Feller} William Feller. \newblock {\em An introduction to probability theory and its applications.  {V}ol. {I}}. \newblock Third edition. John Wiley \& Sons, Inc., New York-London-Sydney,  1968. \MR{0228020}

\bibitem{Fitz}
Patrick Fitzsimmons.
\newblock Excursions above the minimum for diffusions.
\newblock 1985.
\newblock \ARXIV{1308.5189}.

\bibitem{FPY} Patrick Fitzsimmons, Jim Pitman, and Marc Yor. \newblock Markovian bridges: construction, {P}alm interpretation, and splicing. \newblock In {\em Seminar on {S}tochastic {P}rocesses, 1992 ({S}eattle, {WA},  1992)}, volume~33 of {\em Progr. Probab.}, pages 101--134. Birkh\"auser  Boston, Boston, MA, 1993. \MR{1278079}

\bibitem{Forman} Noah Forman. \newblock {\em Instruction sets for walks and the quantile path  transformation}. \newblock ProQuest LLC, Ann Arbor, MI, 2013. \newblock Thesis (Ph.D.)--University of California, Berkeley. \MR{3232248}

\bibitem{Fourati} Sonia Fourati. \newblock Vervaat et {L}\'evy. \newblock {\em Ann. Inst. H. Poincar\'e Probab. Statist.}, 41(3):461--478,  2005. \MR{2139029}

\bibitem{GR} I.~S. Gradshteyn and I.~M. Ryzhik. \newblock {\em Table of integrals, series, and products}. \newblock Elsevier/Academic Press, Amsterdam, seventh edition, 2007. \newblock Translated from the Russian, Translation edited and with a preface by  Alan Jeffrey and Daniel Zwillinger, With one CD-ROM (Windows, Macintosh and  UNIX). \MR{2360010}

\bibitem{Iglehart} Donald~L. Iglehart. \newblock Random walks with negative drift conditioned to stay positive. \newblock {\em J. Appl. Probability}, 11:742--751, 1974. \MR{0368168}

\bibitem{Imhof} J.-P. Imhof. \newblock Density factorizations for {B}rownian motion, meander and the  three-dimensional {B}essel process, and applications. \newblock {\em J. Appl. Probab.}, 21(3):500--510, 1984. \MR{0752015}

\bibitem{LeGall} Jean-Fran{\c{c}}ois Le~Gall and Mathilde Weill. \newblock Conditioned {B}rownian trees. \newblock {\em Ann. Inst. H. Poincar\'e Probab. Statist.}, 42(4):455--489,  2006. \MR{2242956}

\bibitem{Lupu}
Titus Lupu.
\newblock Poissonian ensembles of loops of one-dimensional diffusions.
\newblock 2013.
\newblock \ARXIV{1302.3773}.

\bibitem{Miermont} Gr{\'e}gory Miermont. \newblock Ordered additive coalescent and fragmentations associated to {L}evy  processes with no positive jumps. \newblock {\em Electron. J. Probab.}, 6:no.\ 14, 33 pp. (electronic), 2001. \MR{1844511}

\bibitem{PPY} Mihael Perman, Jim Pitman, and Marc Yor. \newblock Size-biased sampling of {P}oisson point processes and excursions. \newblock {\em Probab. Theory Related Fields}, 92(1):21--39, 1992. \MR{1156448}

\bibitem{Pitman1999} Jim Pitman. \newblock Brownian motion, bridge, excursion, and meander characterized by  sampling at independent uniform times. \newblock {\em Electron. J. Probab.}, 4:no. 11, 33 pp. (electronic), 1999. \MR{1690315}

\bibitem{Pitman} Jim Pitman. \newblock {\em Combinatorial stochastic processes}, volume 1875 of {\em Lecture  Notes in Mathematics}. \newblock Springer-Verlag, Berlin, 2006. \newblock Lectures from the 32nd Summer School on Probability Theory held in  Saint-Flour, July 7--24, 2002, With a foreword by Jean Picard. \MR{2245368}

\bibitem{PR} Jim Pitman and Nathan Ross. \newblock The greatest convex minorant of {B}rownian motion, meander, and  bridge. \newblock {\em Probab. Theory Related Fields}, 153(3-4):771--807, 2012. \MR{2948693}

\bibitem{PY} Jim Pitman and Marc Yor. \newblock Decomposition at the maximum for excursions and bridges of  one-dimensional diffusions. \newblock In {\em It\^o's stochastic calculus and probability theory}, pages  293--310. Springer, Tokyo, 1996. \MR{1439532}

\bibitem{Port} Sidney~C. Port. \newblock An elementary probability approach to fluctuation theory. \newblock {\em J. Math. Anal. Appl.}, 6:109--151, 1963. \MR{0145592}

\bibitem{Protter} Philip Protter. \newblock {\em Stochastic integration and differential equations}, volume~21 of  {\em Applications of Mathematics (New York)}. \newblock Springer-Verlag, Berlin, second edition, 2004. \newblock Stochastic Modelling and Applied Probability. \MR{2020294}

\bibitem{RY}
Daniel Revuz and Marc Yor.
\newblock {\em Continuous martingales and {B}rownian motion}, volume 293 of
  {\em Grundlehren der Mathematischen Wissenschaften [Fundamental Principles of
  Mathematical Sciences]}.
\newblock Springer-Verlag, Berlin, 1991.

\bibitem{Schweinsberg} Jason Schweinsberg. \newblock Applications of the continuous-time ballot theorem to {B}rownian  motion and related processes. \newblock {\em Stochastic Process. Appl.}, 95(1):151--176, 2001. \MR{1847096}

\bibitem{Uribe2014} Ger{\'o}nimo Uribe~Bravo. \newblock Bridges of {L}\'evy processes conditioned to stay positive. \newblock {\em Bernoulli}, 20(1):190--206, 2014. \MR{3160578}

\bibitem{Vervaat} Wim Vervaat. \newblock A relation between {B}rownian bridge and {B}rownian excursion. \newblock {\em Ann. Probab.}, 7(1):143--149, 1979. \MR{0515820}

\bibitem{Vigon} Vincent Vigon. \newblock ({H}omogeneous) {M}arkovian bridges. \newblock {\em Ann. Inst. Henri Poincar\'e Probab. Stat.}, 47(3):875--916,  2011. \MR{2841078}

\bibitem{Wendel} J.~G. Wendel. \newblock Order statistics of partial sums. \newblock {\em Ann. Math. Statist.}, 31:1034--1044, 1960. \MR{0119225}

\bibitem{Williamsbis} David Williams. \newblock Decomposing the {B}rownian path. \newblock {\em Bull. Amer. Math. Soc.}, 76:871--873, 1970. \MR{0258130}

\bibitem{Williams} David Williams. \newblock Path decomposition and continuity of local time for one-dimensional  diffusions. {I}. \newblock {\em Proc. London Math. Soc. (3)}, 28:738--768, 1974. \MR{0350881}

\bibitem{YY} Ju-Yi Yen and Marc Yor. \newblock {\em Local times and excursion theory for {B}rownian motion}, volume  2088 of {\em Lecture Notes in Mathematics}. \newblock Springer, Cham, 2013. \newblock A tale of Wiener and It{\^o} measures. \MR{3134857}

\bibitem{Yorquant} Marc Yor. \newblock The distribution of {B}rownian quantiles. \newblock {\em J. Appl. Probab.}, 32(2):405--416, 1995. \MR{1334895}

\end{thebibliography}



\ACKNO{The second and the third authors would express their gratitude to Noah Forman for helpful discussions and suggestions throughout the preparation of this work. We thank Patrick Fitzsimmons for his remarks on the path decomposition result. The authors also thank three anonymous referees for their careful reading and pointing out errors in an earlier draft, which leads to several improvements in the current version.}


\end{document}